\documentclass[a4paper,reqno]{amsart}
      \usepackage{etex}
      \usepackage[OT2, T1]{fontenc}
      \usepackage{url}
      \usepackage{amsmath}
      \usepackage{lmodern}
      \usepackage{array}
      \usepackage{graphicx}
      \usepackage{amsfonts}
      \usepackage{amssymb}
      \usepackage{extpfeil}
\usepackage{soul}
\usepackage{xargs}
      \usepackage{xcolor}
      \definecolor{royalblue}{RGB}{64, 106, 212}
      \definecolor{link}{RGB}{11,0,128}
      \definecolor{olivegreen}{RGB}{128, 128, 0}
      \usepackage{aliascnt}
      \usepackage[colorlinks=true, citecolor=royalblue,linkcolor=olivegreen,urlcolor=link]{hyperref}
      \usepackage{amsthm}
      \usepackage{esint}
      \usepackage{epic}
      \usepackage{rotating}
\newcommand*{\isoarrow}[1]{\arrow[#1,"\rotatebox{90}{\(\wt{\quad\;}\)}"
]}
      \usepackage{calligra}
      \usepackage{xkeyval}

      \usepackage{exscale, relsize}
      \usepackage{stackengine}
      \usepackage{scalerel}
      \usepackage{setspace}
      \usepackage{paralist}
      \usepackage{colonequals}		
      \usepackage{mathabx}		
      \usepackage[left]{lineno}
      \usepackage{amscd}
      \usepackage[all,cmtip]{xy}	
      \usepackage{sseq}			
      \usepackage{verbatim}
      \usepackage{parskip}
      \usepackage{microtype}
      \usepackage{ragged2e}
            \usepackage[in]{fullpage}
      \usepackage{graphicx}
      \usepackage{mathrsfs} 			
      \usepackage{bm} 				
      \usepackage{cleveref}
      \usepackage{enumitem}
      \usepackage{moreenum}			
      \usepackage{blindtext}
      \usepackage{tikz-cd}
      \usepackage{tikz}
      \usepackage{subfiles}
      \usetikzlibrary{shapes.geometric}
      \usetikzlibrary{automata, bending,matrix,arrows,decorations.pathmorphing,backgrounds,positioning,fit,petri,arrows.meta}
      \usepackage{textcomp}
      \usepackage{indentfirst}
      \usetikzlibrary{arrows}
      \usepackage{enumitem}
      \tikzset{commutative diagrams/.cd,arrow style=tikz,diagrams={>=latex'}}
      \usepackage{filecontents}
      \usepackage[alphabetic,lite]{amsrefs} 	
      \usepackage{stmaryrd}	
      \usepackage{subfiles}
      \usepackage{ragged2e}
      \usetikzlibrary{arrows}
      \usepackage{stackengine}

\usepackage[colorinlistoftodos,prependcaption,textsize=tiny]{todonotes}
\newcommandx{\unsure}[2][1=]{\todo[linecolor=red,backgroundcolor=red!25,bordercolor=red,#1]{#2}}
\newcommandx{\change}[2][1=]{\todo[linecolor=blue,backgroundcolor=blue!25,bordercolor=blue,#1]{#2}}
\newcommandx{\info}[2][1=]{\todo[linecolor=OliveGreen,backgroundcolor=OliveGreen!25,bordercolor=OliveGreen,#1]{#2}}
\newcommandx{\improvement}[2][1=]{\todo[linecolor=Plum,backgroundcolor=Plum!25,bordercolor=Plum,#1]{#2}}
\newcommandx{\thiswillnotshow}[2][1=]{\todo[disable,#1]{#2}}

\usepackage{marginnote}
\newcommand{\mytodo}[2][]{{%
 \let\marginpar\marginnote
 \reversemarginpar
 \renewcommand{\baselinestretch}{0.8}%
 \todo[#1]{#2}}}

         \newcommand{\GG}{\Gamma}
         
         \newcommand{\GGL}{\Lambda}

         \newcommand{\bA}{\mathbb{A}}

         \newcommand{\bF}{\mathbb{F}}
         \newcommand{\bG}{\mathbb{G}}

         \newcommand{\bP}{\mathbb{P}}

         \newcommand{\bZ}{\mathbb{Z}}

         \newcommand{\cF}{\mathcal{F}}
         
         \newcommand{\cH}{\mathcal{H}}

         \newcommand{\cO}{\mathcal{O}}
         \newcommand{\cP}{\mathcal{P}}
         \newcommand{\cQ}{\mathcal{Q}}

         \newcommand{\cW}{\mathcal{W}}

         \newcommand{\fm}{\mathfrak{m}}
         
         \newcommand{\fp}{\mathfrak{p}}
         \newcommand{\fq}{\mathfrak{q}}

         \newcommand{\sE}{\mathscr{E}}
         \newcommand{\sF}{\mathscr{F}}
         \newcommand{\sG}{\mathscr{G}}
         \newcommand{\sH}{\mathscr{H}}

         \newcommand{\sL}{\mathscr{L}}
         
         \newcommand{\sN}{\mathscr{N}}
         \newcommand{\sO}{\mathscr{O}}

         \newcommand{\ra}{\rightarrow}

         \newcommand{\hra}{\hookrightarrow}

         \newcommand{\wt}{\widetilde}
         \newcommand{\wh}{\widehat}

         \newcommand{\pr}{^{\prime}}
         \newcommand{\prpr}{^{\prime\prime}}

         \newcommand{\ce}{\colonequals}
         \newcommand{\ov}{\overline}
         \newcommand{\un}{\underline}

         \renewcommand{\b}{\textbf}
         \newcommand{\surjects}{\twoheadrightarrow}
         
         \newcommand{\isoto}{\overset{\sim}{\longrightarrow}}
         
         \newcommand{\fppf}{\mathrm{fppf}}		                                           
         \newcommand{\et}{\mathrm{\acute{e}t}}	                                           
         \newcommand{\Zar}{\mathrm{Zar}}		                                               
         \newcommand{\sh}{\mathrm{sh}}		                                              

         \newcommand{\mtg}{\un{\mathrm{Tor}}(G)}
         \newcommand{\ddual}{\vee\!\vee}

         \providecommand{\fps}[1]{[\![#1]\!]}
         
         \providecommand{\SP}[1]{\cite[\href{https://stacks.math.columbia.edu/tag/#1}{#1}]{SP}}

         \providecommand{\f}[2]{\frac{#1}{#2}}
         \providecommand{\fps}[1]{\llbracket#1\rrbracket}

\newextarrow{\xbigtoto}{{15}{15}{15}{12}}
   {\bigRelbar\bigRelbar{\bigtwoarrowsleft\rightarrow\rightarrow}}
         
         \DeclareMathOperator{\coker}{coker}		                       
         \DeclareMathOperator{\Spec}{Spec}		                       
         \DeclareMathOperator{\rad}{rad}			                       
         \DeclareMathOperator{\Hom}{Hom}			                       
         \DeclareMathOperator{\Frac}{Frac}		                       
         \DeclareMathOperator{\SL}{SL}			                       
         \DeclareMathOperator{\Ext}{Ext}		                           	
         
         \DeclareMathOperator{\GL}{GL}		                                                  
         \DeclareMathOperator{\Aut}{Aut}		                                                  
         \DeclareMathOperator{\Lie}{Lie}		                                                  
         \DeclareMathOperator{\Pic}{Pic}		                                                  
         \DeclareMathOperator{\codim}{codim}		                                                  
         \DeclareMathOperator{\Isom}{Isom}	                                             	
         \DeclareMathOperator{\Jac}{Jac}		                                             

         \newcommand{\ba}{\begin{aligned}}
         \newcommand{\ea}{\end{aligned}}
         \newcommand{\be}{\begin{equation}}
         \newcommand{\ee}{\end{equation}}
         \newcommand{\pf}{\begin{proof}}
         \newcommand{\bpf}{\begin{proof}}
         \newcommand{\epf}{\end{proof}}
         \newcommand{\bsol}{\begin{solution}}
         \newcommand{\esol}{\end{solution}}
         \newcommand{\bthm}{\begin{thm}}
         \newcommand{\ethm}{\end{thm}}
         \newcommand{\bthmt}{\begin{thm-tweak}}
         \newcommand{\ethmt}{\end{thm-tweak}}
         \newcommand{\bprop}{\begin{prop}}
         \newcommand{\eprop}{\end{prop}}
         \newcommand{\bcor}{\begin{cor}}
         \newcommand{\ecor}{\end{cor}}
         \newcommand{\bcort}{\begin{cor-tweak}}
         \newcommand{\ecort}{\end{cor-tweak}}
         \newcommand{\brem}{\begin{rem}}
         \newcommand{\erem}{\end{rem}}
         \newcommand{\bremt}{\begin{rem-tweak}}
         \newcommand{\eremt}{\end{rem-tweak}}
         \newcommand{\brems}{\begin{rems} \hfill \begin{enumerate}[label=\b{\thenumberingbase.},ref=\thenumberingbase]}
         
         \newcommand{\erems}{\end{enumerate} \end{rems}}
         \newcommand{\begs}{\begin{egs} \hfill \begin{enumerate}[label=\b{\thenumberingbase.},ref=\thenumberingbase]}
         
         \newcommand{\eegs}{\end{enumerate} \end{egs}}
         \newcommand{\eremstweak}{\end{enumerate} \end{rems-tweak}}
         \newcommand{\eremst}{\end{enumerate} \end{rems-tweak}}
         \newcommand{\blem}{\begin{lemma}}
         \newcommand{\elem}{\end{lemma}}
         \newcommand{\blemt}{\begin{lemma-tweak}}
         \newcommand{\elemt}{\end{lemma-tweak}}
         \newcommand{\bconj}{\begin{conj}}
         \newcommand{\econj}{\end{conj}}
         \newcommand{\bprob}{\begin{Problem}}
         \newcommand{\eprob}{\end{Problem}}
         \newcommand{\bpropt}{\begin{prop-tweak}}
         \newcommand{\epropt}{\end{prop-tweak}}
         \newcommand{\bq}{\begin{Q}}
         \newcommand{\eq}{\end{Q}}
         \newcommand{\benum}{\begin{enumerate}[label={{\upshape(\alph*)}}]}
         \newcommand{\benuma}{\begin{enumerate}[label={{\upshape(\arabic*)}}]}
         \newcommand{\benumr}{\begin{enumerate}[label={{\upshape(\roman*)}}]}
         \newcommand{\eenum}{\end{enumerate}}
         \newcommand{\bc}{}
         \newcommand{\bd}{\begin{defn}}
         \newcommand{\ed}{\end{defn}}
         \newcommand{\bque}{\begin{que}}
         \newcommand{\eque}{\end{que}}
         \newcommand{\bfct}{\begin{fact}}
         \newcommand{\efct}{\end{fact}}

         \newcommand{\beg}{\begin{eg}}
         \newcommand{\eeg}{\end{eg}}
         \newcommand{\begt}{\begin{eg-tweak}}
         \newcommand{\eegt}{\end{eg-tweak}}
         \newcommand{\bcl}{\begin{claim}}
         \newcommand{\ecl}{\end{claim}}
         \newcommand{\bclt}{\begin{cl-tweak}}
         \newcommand{\eclt}{\end{cl-tweak}}

         \newcommand{\x}{\text}

         \newcommand{\q}{\quad}
         \newcommand{\qq}{\quad\quad}
         \newcommand{\qqq}{\quad\quad\quad}

         \newcommand{\tst}{\textstyle}

         \newcommand{\sHom}{\mathscr{H}\! om}
         
         \newcommand{\bconjt}{\begin{conj-tweak}}
         \newcommand{\econjt}{\end{conj-tweak}}

         \renewcommand{\Isom}{\un{\mathrm{Isom}}}

\providecommand{\SPD}[2]{\cite[\href{https://stacks.math.columbia.edu/tag/#1}{#1}, \href{https://stacks.math.columbia.edu/tag/#2}{#2}]{SP}}

\newcommand*{\QED}{\hfill\ensuremath{\qed}}

\tikzset{
    labl/.style={anchor=south, rotate=90, inner sep=.5mm}
}

\newaliascnt{numberingbase}{subsection}
\numberwithin{equation}{numberingbase}

\newcommand{\statementskip}{9pt plus 2pt minus 1pt}
\newcommand{\subsectionbeforeskip}{9pt plus 2pt minus 1pt}
\newcommand{\subsectionafterskip}{4pt plus 1pt minus .5pt}
\makeatletter
\def\subsection{\@startsection{subsection}{2}%
  \z@{\subsectionbeforeskip}{\subsectionafterskip}%
  {\normalfont\bfseries}}
\def\subsubsection{\@startsection{subsubsection}{3}%
  \z@{7pt plus 1.5pt minus 1pt}{3pt plus 1pt minus .5pt}%
  {\normalfont\itshape}}
\makeatother

\newtheoremstyle{thms}{\statementskip}{\statementskip}{\itshape}{}{\bfseries}{.}{ }{}
\theoremstyle{thms}
\newaliascnt{conj}{numberingbase}
\newtheorem{conj}[conj]{Conjecture}
\aliascntresetthe{conj}
\newaliascnt{corollary}{numberingbase}

\aliascntresetthe{corollary}
\newaliascnt{cor}{numberingbase}
\newtheorem{cor}[cor]{Corollary}
\aliascntresetthe{cor}
\newaliascnt{lemma}{numberingbase}
\newtheorem{lemma}[lemma]{Lemma}
\aliascntresetthe{lemma}
\newaliascnt{sublemma}{equation}

\aliascntresetthe{sublemma}
\newaliascnt{prop}{numberingbase}
\newtheorem{prop}[prop]{Proposition}
\aliascntresetthe{prop}
\newaliascnt{proposition}{numberingbase}

\aliascntresetthe{proposition}
\newaliascnt{Q}{numberingbase}
\newtheorem{Q}[Q]{Question}
\aliascntresetthe{Q}
\newaliascnt{thm}{numberingbase}
\newtheorem{thm}[thm]{Theorem}
\aliascntresetthe{thm}
\newaliascnt{theorem}{numberingbase}

\aliascntresetthe{theorem}
\newaliascnt{variant}{numberingbase}
\newtheorem{variant}[variant]{Variant}
\aliascntresetthe{variant}

\newtheoremstyle{claims}{\statementskip}{\statementskip}{}{}{\itshape}{.}{ }{}
\theoremstyle{claims}
\newaliascnt{claim}{equation}
\newtheorem{claim}[claim]{Claim}
\aliascntresetthe{claim}
\newaliascnt{cl-tweak}{subsubsection}
\newtheorem{cl-tweak}[cl-tweak]{Claim}
\aliascntresetthe{cl-tweak}
\Crefname{cl-tweak}{Claim}{Claims}

\newtheoremstyle{defs}{\statementskip}{\statementskip}{}{}{\bfseries}{.}{ }{}
\theoremstyle{defs}
\newaliascnt{defn}{numberingbase}
\newtheorem{defn}[defn]{Definition}
\aliascntresetthe{defn}
\newaliascnt{definition}{numberingbase}

\aliascntresetthe{definition}
\newaliascnt{eg}{numberingbase}
\newtheorem{eg}[eg]{Example}
\aliascntresetthe{eg}
\newaliascnt{example}{numberingbase}

\aliascntresetthe{example}
\newtheorem*{egs}{Examples}
\newaliascnt{rem}{numberingbase}
\newtheorem{rem}[rem]{Remark}
\aliascntresetthe{rem}
\newaliascnt{remark}{numberingbase}

\aliascntresetthe{remark}
\newtheorem*{rems}{Remarks}

\Crefname{claim}{Claim}{Claims}
\Crefname{sublemma}{Lemma}{Lemmas}
\Crefname{conj}{Conjecture}{Conjectures}
\Crefname{cor}{Corollary}{Corollaries}
\Crefname{defn}{Definition}{Definitions}
\Crefname{eg}{Example}{Examples}
\Crefname{prop}{Proposition}{Propositions}
\Crefname{Q}{Question}{Questions}
\Crefname{rem}{Remark}{Remarks}
\Crefname{thm}{Theorem}{Theorems}
\Crefname{variant}{Variant}{Variants}
\crefname{claim}{claim}{claims}
\crefname{sublemma}{lemma}{lemmas}
\crefname{conj}{conjecture}{conjectures}
\crefname{cor}{corollary}{corollaries}
\crefname{corollary}{corollary}{corollaries}
\crefname{defn}{definition}{definitions}
\crefname{definition}{definition}{definitions}
\crefname{eg}{example}{examples}
\crefname{example}{example}{examples}
\crefname{lemma}{lemma}{lemmas}
\crefname{prop}{proposition}{propositions}
\crefname{proposition}{proposition}{propositions}
\crefname{Q}{question}{questions}
\crefname{rem}{remark}{remarks}
\crefname{remark}{remark}{remarks}
\crefname{thm}{theorem}{theorems}
\crefname{theorem}{theorem}{theorems}
\crefname{variant}{variant}{variants}
\Crefname{corollary}{Corollary}{Corollaries}
\Crefname{definition}{Definition}{Definitions}
\Crefname{example}{Example}{Examples}
\Crefname{lemma}{Lemma}{Lemmas}
\Crefname{proposition}{Proposition}{Propositions}
\Crefname{remark}{Remark}{Remarks}
\Crefname{theorem}{Theorem}{Theorems}

\theoremstyle{thms}
\newaliascnt{thm-tweak}{subsection}
\newtheorem{thm-tweak}[thm-tweak]{Theorem}
\aliascntresetthe{thm-tweak}
\crefname{thm-tweak}{theorem}{theorems}
\Crefname{thm-tweak}{Theorem}{Theorems}
\newaliascnt{lemma-tweak}{subsection}
\newtheorem{lemma-tweak}[lemma-tweak]{Lemma}
\aliascntresetthe{lemma-tweak}
\crefname{lemma-tweak}{lemma}{lemmas}
\Crefname{lemma-tweak}{Lemma}{Lemmas}
\newaliascnt{cor-tweak}{subsection}
\newtheorem{cor-tweak}[cor-tweak]{Corollary}
\aliascntresetthe{cor-tweak}
\crefname{cor-tweak}{corollary}{corollaries}
\Crefname{cor-tweak}{Corollary}{Corollaries}
\newaliascnt{prop-tweak}{subsection}
\newtheorem{prop-tweak}[prop-tweak]{Proposition}
\aliascntresetthe{prop-tweak}
\crefname{prop-tweak}{proposition}{propositions}
\Crefname{prop-tweak}{Proposition}{Propositions}
\newaliascnt{conj-tweak}{subsection}
\newtheorem{conj-tweak}[conj-tweak]{Conjecture}
\aliascntresetthe{conj-tweak}
\crefname{conj-tweak}{conjecture}{conjectures}
\Crefname{conj-tweak}{Conjecture}{Conjectures}

\theoremstyle{defs}
\newaliascnt{defn-tweak}{subsection}
\newtheorem{defn-tweak}[defn-tweak]{Definition}
\aliascntresetthe{defn-tweak}
\crefname{defn-tweak}{definition}{definitions}
\Crefname{defn-tweak}{Definition}{Definitions}
\newaliascnt{eg-tweak}{subsection}
\newtheorem{eg-tweak}[eg-tweak]{Example}
\aliascntresetthe{eg-tweak}
\crefname{eg-tweak}{example}{examples}
\Crefname{eg-tweak}{Example}{Examples}
\newtheorem*{rems-tweak}{Remarks}
\newaliascnt{rem-tweak}{subsection}
\newtheorem{rem-tweak}[rem-tweak]{Remark}
\aliascntresetthe{rem-tweak}
\crefname{rem-tweak}{remark}{remarks}
\Crefname{rem-tweak}{Remark}{Remarks}

\newtheoremstyle{subsection-tweak}
   {\subsectionbeforeskip}
   {\subsectionafterskip}%
   {}
   {}%
   {\bfseries}
   {}%
   {\newline}
   {\thmnumber{\@{#1}{}\@{#2}.}%
    \thmnote{~{\bfseries#3.}}}

\makeatletter
\newcommand{\ppafterheading}{%
  \@ifnextchar\label{\ppafterheadinglabel}{\ignorespaces}%
}
\def\ppafterheadinglabel\label#1{\label{#1}\ignorespaces}
\makeatother

\newenvironment{pp}[1][]{%
  \par\addvspace{\subsectionbeforeskip}%
  \refstepcounter{subsection}%
  \noindent\textbf{\thesubsection.\if\relax\detokenize{#1}\relax\else\ #1.\fi}%
  \par\nobreak\vspace{\subsectionafterskip}\noindent\ppafterheading
}{%
  \par\addvspace{\subsectionafterskip}%
}
\newcommand{\bpp}{\begin{pp}}
\newcommand{\epp}{\end{pp}}
\newaliascnt{pp-t}{subsubsection}
\newtheorem{pp-t}[pp-t]{}
\aliascntresetthe{pp-t}
\crefname{pp-t}{section}{sections}
\Crefname{pp-t}{Section}{Sections}

\theoremstyle{subsection-tweak}
\newaliascnt{conventions}{subsection}

\aliascntresetthe{conventions}
\crefname{conventions}{section}{sections}
\Crefname{conventions}{Section}{Sections}

\theoremstyle{subsection-tweak}
\newtheorem{pp-tweak}{}
\crefname{pp-tweak}{section}{sections}
\Crefname{pp-tweak}{Section}{Sections}

      \makeatletter
      \def\@tocline#1#2#3#4#5#6#7{
          \begingroup
          \@ifempty{#4}{}{}

          \parindent\z@ \leftskip#3\relax \advance\leftskip\@tempdima\relax
          #5\hskip-\@tempdima
            \ifcase #1
             \or\or \hskip 2em \or \hskip 1em \else \hskip 3em \fi%
            #6\nobreak\relax
          \dotfill\hbox to\@pnumwidth{\@tocpagenum{#7}}\par
          \nobreak
          \endgroup
        }
       \def\l@section{\@tocline{1}{0pt}{1pc}{}{}}

      \renewcommand{\tocsection}[3]{%
        \indentlabel{\@ifnotempty{#2}{\makebox[1.3em][l]{%
          \ignorespaces#1 \bfseries{#2}.\hfill}}}\bfseries{#3}
          \vspace{-3.5pt}}

      \renewcommand{\tocsubsection}[3]{%
        \indentlabel{\@ifnotempty{#2}{\hspace*{-0.5em}\makebox[2.1em][l]{%
          \ignorespaces#1#2.\hfill}}}#3
          \vspace{-4.5pt}}

\makeatother

\makeatletter 
\newcommand\appendix@section[1]{%
  \refstepcounter{section}%
  \orig@section*{Appendix \@Alph\c@section. #1}%
}
\let\orig@section\section
\g@addto@macro\appendix{\let\section\appendix@section}
\makeatother

      \DeclareMathVersion{normal2}

\author{Ning Guo}
\address{Institute for Advanced Study in Mathematics of Harbin Institute of Technology, Harbin, China}
\email{ninguo@hit.edu.cn}
\author{Fei Liu}
\address{Department of Mathematics, Southern University of Science and Technology, Shenzhen, China}
\email{liufei54@pku.edu.cn}
\date{\today}

\def\UTFviii@defined#1{%
  \ifx#1\relax
      \PackageError{inputenc}{Unicode\space char\space\expandafter
                              \UTFviii@splitcsname\string#1\relax
                              \MessageBreak
                              not\space set\space up\space
                              for\space use\space with\space LaTeX}\@eha
  \else\expandafter
    #1%
  \fi
}

\makeatletter
\def\UTFviii@defined#1{%
  \ifx#1\relax
      ?%
  \else\expandafter
    #1%
  \fi
}

\makeatother


\begin{document}

\title{Grothendieck--Serre for constant reductive group schemes}
\begin{abstract}
The Grothendieck--Serre conjecture predicts that on a regular local ring there is no nontrivial torsor under a reductive group scheme that becomes trivial over the fraction field.
While this conjecture has been proven in the equicharacteristic case, it remains open in the mixed characteristic case.
In this article, we establish a generalised version of the conjecture over Pr\"ufer bases for constant reductive group schemes.
In particular, the Noetherian case of our main result settles the constant, unramified case of the Grothendieck--Serre conjecture.
Along the way, inspired by the recent article by $\check{\mathrm{C}}$esnavi$\check{\mathrm{c}}$ius \cite{Ces24}, we also prove several versions of the Nisnevich conjecture in our context.
\end{abstract}
\maketitle

\vspace{-25pt}



\hypersetup{
    linktoc=page,     
}



\section{Introduction} 
The Grothendieck--Serre conjecture predicts that a torsor under a reductive group scheme $G$ over a semilocal regular ring $A$ is trivial if it trivializes over the total ring of fractions.
Equivalently, the map
\[
   H^1_{\et}(A,G)\ra H^1_{\et}(\Frac A,G)
\]
between nonabelian cohomology pointed sets has trivial kernel.
The conjecture has been proven when $\dim A\leq 1$, or when $A$ contains a field (i.e., when $A$ is of equicharacteristic).
When $A$ is of mixed characteristic, the conjecture is widely open, though several unramified cases have been settled. 
For a detailed review of the history, see the survey \cite{Pan18}*{\S 5} and the summary \S\ref{history} below, which includes recent developments.
The first goal of this article is to establish a new mixed characteristic case of the Grothendieck--Serre conjecture in a much more general non-Noetherian setup.

\bthmt [Theorem~\ref{G-S for constant reductive gps}~\ref{G-S for constant reductive gps i}]\label{main thm}
For a Pr\"ufer ring $R$, an irreducible, affine, smooth $R$-scheme $X$, and a reductive $R$-group scheme $G$, every generically trivial $G$-torsor over $X$ is Zariski semilocally trivial. In particular, the sequence
\[
   \text{$1\ra H^1_{\Zar}(X,G)\ra H^1_{\et}(X,G)\ra H^1_{\et}(K(X),G)$\q is exact,}
\]
where $K(X)$ is the function field of $X$.
In other words, for every semilocal ring $A$ of $X$, we have
\[
  \ker\,(H^1_{\et}(A,G)\ra H^1_{\et}(\Frac A, G))=\{\ast\}.
\]
\ethmt

We recall that a ring is \emph{Pr\"ufer} if all of its local rings are valuation rings, that is, normal domains whose ideals are totally ordered by inclusion.
Thus, the Noetherian Pr\"ufer rings are exactly the Dedekind rings.
Pr\"ufer rings are ubiquitous in nonarchimedean analytic geometry and perfectoid theory.
For instance, all perfectoid fields are \emph{non-discretely} valued, and their integer rings are \emph{non-Noetherian} valuation rings.
Moreover, in the recently introduced v- and arc-topologies, covers are detected by maps from spectra of valuation rings. This provides a powerful way to reduce many questions to their analogues over valuation rings. These developments motivate the study of geometry over Pr\"ufer rings and, in particular, of the Grothendieck--Serre problem in this setting.

The novelty and significance of  \Cref{main thm} and its approach are highlighted in the following four key aspects.

$\bullet$ Even when $X$ has relative dimension 0 over $R$ (\emph{i.e.}, when $A$ is semilocal Pr\"uferian) and $G$ is merely a torus, as pointed out by Colliot-Thélène, it is unfeasible to apply the argument in \cite{CTS78} anymore, even though it works in the case of discrete valuation rings.
The obstruction to that argument is the failure of $\Ext$ functors to commute with limits.
To circumvent this, we employ the coniveau spectral sequence, which not only reduces the problem to the resolved valuation ring case but also extends the result to higher relative dimensional cases by using purity for torsors under tori on smooth schemes over Prüfer bases when $G$ is a torus.
Note that \Cref{main thm} recovers the main result, \cite{Guo24}, of the first author, where $A$ is an arbitrary valuation ring. 

$\bullet$ We prove a new geometric presentation lemma in the mixed characteristic setting, a result that should also be useful in other contexts, such as in $K$-theory and $\bA^1$-homotopy theory (see also \cite{GL25}*{\S1.14}). 
Specifically, if $X$ is smooth of relative dimension $d>0$ over a discrete valuation ring $R$ (or, more generally, a valuation ring $R$) and $Y\subset X$ is a fibrewise nowhere dense closed subscheme, then, Zariski-locally on $X$, there exists
an \'etale map $(X,Y)\to (X\pr,Y\pr)$ to a similar pair $(X\pr,Y\pr)$, with $Y\xrightarrow{\sim} X\times_{X\pr}Y\pr\xrightarrow{\sim} Y\pr$ and $X\pr\subset \bA^d_{R}$ an open subset. For the precise statement, see \Cref{variant of Lindel's lem}.

$\bullet$ The Noetherian case of \Cref{main thm} already resolves a new case of the Grothendieck--Serre conjecture. Namely, the conjecture holds when $A$ is a regular semilocal ring essentially smooth over a discrete valuation ring and $G$ descends to a reductive group scheme over that ring,  
This special case is one of the state-of-the-art results on the Grothendieck--Serre conjecture.
The approach leverages the geometric result mentioned earlier to reduce the problem to the case where $X\subset \mathbb{A}_R^d$ is an open subset, a scenario in which $X$ admits a \emph{smooth projective} compactification $\ov{X}$.
The constancy of the reductive group $G$ is crucial for this reduction.
Once such a compactification $\ov{X}$ exists, $G$ need not be constant; it is enough for $G$ to be defined over $\ov X$, as shown in the following variant.
\bthmt[Theorem~\ref{torsors-Sm proj base}\ref{loc-gen-trivial-sm-proj}]\label{Main-Sm proj base}
   For a semilocal Pr\"{u}fer domain $R$, an irreducible, smooth, projective $R$-scheme $X$, a reductive $X$-group scheme $G$, and a semilocal ring $A$ of $X$, the natural map
   \[
     \text{$H^1_{\et}(A,G)\ra H^1_{\et}(\Frac A, G)$ \q has trivial kernel.}
   \]
\ethmt
The special case where $R$ is a discrete valuation ring, $G$ descends to a reductive $R$-group scheme, and the generically trivial torsor is defined on all of $X$ was later obtained independently by the first author and Panin in a subsequent preprint \cite{GP23b}, as a consequence of \cite{GP23a}, \cite{GP25}, and \cite{PS23}. In contrast, our proof of \Cref{Main-Sm proj base} is shorter and more direct, relying only on the presentation lemma \cite{Ces22a}*{Variant~3.7} and the purity result \cite{CTS79}*{théorème~6.13} of Colliot-Th\'el\`ene and Sansuc.

$\bullet$ Moreover, \Cref{main thm} implies a stronger assertion in the Noetherian case: when $A$ is a semilocal ring that is geometrically regular over some Dedekind ring $R$, and $G$ is a reductive $R$-group scheme, the Grothendieck--Serre conjecture still holds.
This follows from the Popescu's theorem \SP{07GC}, which allows us to use limit formalism to reduce to the case resolved in \Cref{main thm}.

\bpp[Known cases of the Grothendieck--Serre conjecture]
\label{history}
Since proposed by Serre \cite{Ser58}*{page~31} and Grothendieck \cite{Gro58}*{pages 26--27, remarques~3}, \cite{Gro68}*{remarques~1.11~a)}, the Grothendieck--Serre conjecture has been studied extensively, and several cases are known, as summarized below.
\benumr
\item The case when $G$ is a torus was proved by Colliot-Th\'el\`ene and Sansuc in \cite{CTS87}.
\item The case when $A$ is a discrete valuation ring was addressed by Nisnevich in \cite{Nis82} and \cite{Nis84}. This result was improved and extended to the semilocal Dedekind case in \cite{Guo22}.
Several special cases were proved in \cite{Har67}, \cite{BB70}, \cite{BrT3} for discrete valuation rings, and in \cite{PS16}, \cite{BVG14}, \cite{BFF17}, \cite{BFFH20} for the semilocal Dedekind case.
\item The case when $A$ is Henselian local was settled in \cite{BB70} and \cite{CTS79}*{assertion~6.6.1}. In this case, Hensel's lemma reduces the triviality of $G$-torsors to that over the residue fields. By inducting on $\dim A$, one can eventually reach the resolved case of Nisnevich.
\item The equicharacteristic case, that is, when $A$ contains a field $k$, was established by Fedorov and Panin \cite{FP15} when $k$ is infinite (see also \cite{PSV15}, \cite{Pan20b} for crucial techniques) and by Panin \cite{Pan20a} when $k$ is finite. These arguments are streamlined and simplified in \cite{Fed22a}. For constant groups over fields, Bouthier--\v{C}esnavi\v{c}ius--Scavia \cite{BCS25} recently proved the stronger statement that, for a smooth $k$-group scheme $G$ and a smooth $k$-variety $X$, every generically trivial $G$-torsor on $X$ is Zariski semilocally trivial.
Prior to these works, several equicharacteristic subcases were proven in \cite{Oja80}, \cite{CTO92}, \cite{Rag94}, \cite{PS97}, \cite{Zai00}, \cite{OP01}, \cite{OPZ04}, \cite{Pan05}, \cite{Zai05}, \cite{Che10}, \cite{PSV15}.
\item In mixed characteristic, \v{C}esnavi\v{c}ius \cite{Ces22a}*{Theorem~5.3.1} settled the quasi-split case for semilocal regular algebras over a Dedekind ring whose fibres are geometrically regular and whose residue field extensions at maximal ideals satisfy the separability condition in that theorem. This contains the unramified regular local case.
 Fedorov \cite{Fed22b} previously proved the split case under additional assumptions.
\v{C}esnavi\v{c}ius \cite{Ces24}*{Theorem~1.3} settled a generalised Nisnevich conjecture under certain conditions, which specializes to the Grothendieck--Serre conjecture in the equal and mixed characteristic cases that were previously proven  in \cite{FP15}, \cite{Pan20a}, \cite{Ces22a}.
\v{C}esnavi\v{c}ius--Fedorov \cite{CF23} (after Fedorov's result \cite{Fed23}) proved the case when $A$ is unramified and $G^{\mathrm{ad}}$ is totally isotropic.
Their proof replaces affine Grassmannians by the geometry of $\mathrm{Bun}_G$ and develops reembedding and excision techniques for relative curves in which finiteness is weakened to quasi-finiteness.

\item There are sporadic cases where either $A$ or $G$ are special, possibly with additional mixed characteristic conditions, which have been settled in the literature. These include references such as \cite{Gro68}*{Remarque~1.11~a)}, \cite{Oja82}, \cite{Nis89}, \cite{Fed22b}, \cite{Fir23}, \cite{BFFP22}, \cite{Pan19b}.
\eenum
In proving \Cref{main thm}, we use only our Pr\"ufer analogues of the toral case treated in \cite{CTS87} (see \Cref{G-S type results for mult type}) and the semilocal Dedekind case settled in \cite{Guo22} (see \Cref{G-S over semi-local prufer}); no other known case of the Grothendieck--Serre conjecture is needed.
Our proof is outlined in \S\ref{intro-outline-pf of main thm}.
\epp

\bpp[Nisnevich's purity conjecture]
\label{intro-Nisnevich purity}
The second topic of this article is Nisnevich's purity, which requires a total isotropicity assumption on the adjoint quotient of the group schemes involved.

A semisimple group scheme $G$ over $S$ is \emph{totally isotropic} at $s\in S$ if every factor $G_i$ in the decomposition \cite{SGA3IIInew}*{exposé~XXIV, proposition~5.10~(i)}
\[
   \tst G^{\mathrm{ad}}_{\sO_{S,s}}\cong \prod_i\mathrm{Res}_{R_i/\sO_{S,s}}(G_i)
\]
contains a copy of $\bG_{m,R_i}$, where each $\sO_{S,s}\to R_i$ is a finite \'etale cover, and $G_i$ is a simple\footnote{This means that the geometric fibres are reductive algebraic groups having connected Dynkin diagrams.} adjoint $R_i$-group scheme (of constant type); equivalently, every $G_i$ has a parabolic $R_i$-subgroup that is fibrewise proper, see \cite{SGA3IIInew}*{expos\'e XXVI, corollaire 6.12}.
If this holds for all $s\in S$, then $G$ is \emph{totally isotropic}. 
For instance, all quasi-split (in particular, split) semisimple group schemes are totally isotropic.
\epp

The following conjecture dates back to Nisnevich \cite{Nis89}*{Conjecture~1.3} and was later modified in light of Fedorov's anisotropic counterexamples \cite{Fed22b}*{Proposition~4.1}.

\bconjt[Nisnevich]
\label{Nisnevich-introd}
For a regular semilocal ring $A$, an element $r\in A$ such that $r\notin \fm^2$ for every maximal ideal $\fm\subset A$, and a reductive $A$-group scheme $G$ such that $G_{A/rA}^{\mathrm{ad}}$ is totally isotropic, every generically trivial $G$-torsor on $A[\f{1}{r}]$ is trivial, that is, the natural map
\[
     \text{$H^1_{\et}(A[1/r],G)\ra H^1_{\et}(\Frac A, G)$ \q has trivial kernel.}
   \]
\econjt

Note that this conjecture is already nontrivial for $(A,\fm)$ a regular local ring and $G=\GL_n$, a case known as the Quillen conjecture: if $r\notin \fm^2$, then every finite projective $A[\frac{1}{r}]$-module ought to be free.
In fact, its special case when $A$ is a local ring of a regular affine variety over a field and $G=\GL_n$ was settled by Bhatwadekar--Rao in \cite{BR83} and was subsequently extended to arbitrary regular local rings containing fields by Popescu \cite{Pop02}*{Theorem~1}.
 Nisnevich in \cite{Nis89} proved the conjecture in dimension $2$, assuming that $A$ is a local ring with an infinite residue field and that $G$ is quasi-split.
For the state of the art, the conjecture was settled in the equicharacteristic case and in several mixed characteristic cases by {\v{C}}esnavi{\v{c}}ius in \cite{Ces24}*{Theorem~1.3}. 
Previously, Fedorov \cite{Fed24} proved the case when $A$ contains an infinite field.
The toral case and some low dimensional cases are known and surveyed in \cite{Ces22}*{\S\S3.4.2(1)}, including Gabber's result \cite{Gab81}*{Chapter~I, Theorem~1} for the local case $\dim A\leq 3$ when $G$ is either $\GL_n$ or $\mathrm{PGL}_n$.

In this article, we prove several variants of \Cref{Nisnevich-introd}, over Pr\"ufer bases, such as the following.

\bthmt
 Let $R$ be a semilocal Pr\"{u}fer domain, let $X$ be an irreducible, smooth, quasi-projective $R$-scheme, and let $G$ be a reductive $X$-group scheme. Let $A$ be a semilocal ring of $X$.
Let $r\in R\backslash \{0\}$. Suppose that either
   \benumr
   \item \emph{(\Cref{torsors-Sm proj base}~\ref{Nis-sm-proj})} $X$ is $R$-projective and $G_{A[\f{1}{r}]}^{\mathrm{ad}}$ is totally isotropic; or 
   \item \emph{(\Cref{G-S for constant reductive gps}~\ref{G-S for constant reductive gps ii})} $G$ descends to a reductive $R$-group $G\pr$ with $(G\pr_{R[\f{1}{r}]})^{\mathrm{ad}}$ totally isotropic.
   \eenum
 Then every generically trivial $G$-torsor over $A[\f{1}{r}]$ is trivial, that is, we have
   \[
   \tst \ker\,(H^1_{\et}(A[\f{1}{r}],G)\ra H^1_{\et}(\Frac A,G))=\{\ast\}.
   \]
\ethmt

\bpp[Outline of the proof of \Cref{main thm}]
\label{intro-outline-pf of main thm}
In the setting of \Cref{main thm}, following \v{C}esnavi\v{c}ius and earlier ideas of Fedorov--Panin, one tries to use a Gabber--Quillen type geometric presentation lemma to fibre a suitable open neighbourhood $U\subset X$ of a finite subset $\mathbf{x}\subset X$ into smooth affine curves $U\to S$ over an open
\[
S\subset \bA_{R}^{\dim(X/R)-1}
\]
in such a way that a given small\footnote{Typically, it is not so small, as it is only $R$-fibrewise of codimension $\ge 1$ in $X$.} closed subscheme $Y\subset X$ restricts to an $S$-\emph{finite} (not merely $S$-quasifinite) closed subscheme of $U$. For us, $Y$ is a closed subscheme away from which the generically trivial torsor in question is trivial and which contains no irreducible component of a special fibre of $X$. Unfortunately, the $S$-finiteness mentioned above could be very difficult (and even impossible) to achieve in general, especially when working over mixed characteristic rings.  In view of the mixed characteristic variant of the Gabber--Quillen type geometric presentation lemma \cite{Ces22a}*{Variant 3.7}, one way to overcome this difficulty is to compactify $X$ to a projective, flat scheme $\overline{X}$ over $R$ such that the boundary $\overline{Y}\backslash Y$ of $Y$ has $R$-fibrewise codimension $\ge 2$ in $\overline{X}$.

Assume that we can find such a compactification $\overline{X}$ which is even $R$-smooth. Although this assumption may appear restrictive, it suffices for our later application where $\overline{X} = \mathbb{P}_{R}^d$. A key consequence of the purity for torsors under reductive group schemes on schemes smooth over Pr\"ufer schemes (\Cref{extends across codim-2 points}) is that we can enlarge the domain $X$ of our torsor such that $\overline{X}\backslash X$ is $R$-fibrewise of codimension $\ge 2$ in $\overline{X}$. A more precise statement can be found in \Cref{extend generically trivial torsors}. With this enlargement, we can ensure that $\overline{Y}\backslash Y$ is also $R$-fibrewise of codimension $\ge 2$ in $\overline{X}$, making the geometric presentation lemma mentioned above applicable. Consequently, following the arguments of Fedorov--Panin and \v{C}esnavi\v{c}ius, one can reduce, via excision and patching techniques, to studying the triviality of torsors on the relative affine line---a scenario analysed in detail in \Cref{triviality on sm rel. affine curves}.

The original smooth $R$-scheme $X$, however, need not admit a smooth projective compactification. To address this, we prove the following result, which is interesting in its own right: if $Y$ is a fibrewise nowhere dense closed subset of $X$, then, Zariski-locally on $X$, the pair $(Y,X)$ can be presented as an elementary étale neighbourhood of a similar pair $(Y\pr,X\pr)$, where $X\pr$ is an open subset of some projective $R$-space; see \Cref{variant of Lindel's lem} for a more precise statement. 
By using standard gluing techniques, we can replace $(Y,X)$ with $(Y\pr,X\pr)$ and study generically trivial torsors on $X\pr$. 
Moreover, since now $X\pr$ has smooth projective compactifications, namely, the projective $R$-spaces, we can reduce to the situation already settled in the previous paragraph.
\epp

\bpp[Notations and conventions]
Throughout this paper, we work with commutative rings with units. For a
point $s$ of a scheme (resp., for a prime ideal $\mathfrak{p}$ of a ring), we let $k_s$ (resp., $k_{\mathfrak{p}}$) denote its residue
field. Given a global section $a$ of a scheme $S$, we use $S[\frac{1}{a}]$ to denote the open subset of $S$ where $a$ does not vanish. For a ring $A$, we let $\Frac A$ denote its total ring of fractions.
For a morphism of schemes $S\pr\to S$, we denote the base change functor from $S$ to $S\pr$ by $(-)_{S\pr}$. If $S\pr=\Spec R\pr$ is affine, we also write $(-)_{R\pr }$ for $(-)_{S\pr}$.

 Let $S$ be a scheme, and let $G$ be an $S$-group scheme. Given an $S$-scheme $T$, by a $G$-torsor over $T$ we shall mean a $G_T\ce G\times_ST$-torsor. We denote by $\textbf{Tors}(S_{\fppf},G)$ (resp., $\textbf{Tors}(S_{\et},G)$) the groupoid of $G$-torsors on $S$ that are fppf-locally (resp., \'etale-locally) trivial. Specifically, if $G$ is $S$-smooth (e.g., if $G$ is $S$-reductive), then every fppf-locally trivial $G$-torsor is \'etale-locally trivial, so we have
$$\textbf{Tors}(S_{\fppf},G)= \textbf{Tors}(S_{\et},G).$$

Let $X$ be a scheme. We write $(-)^{\vee}\ce \sHom_{\sO_X}(-,\sO_X)$ for duals of $\sO_X$-modules. A coherent $\sO_X$-module $\sF$ is \emph{reflexive} if the canonical map $\sF\ra \sF^{\ddual}$ is an isomorphism. The Picard groupoid of $X$, denoted by $\mathbf{Pic}_X$, is the groupoid of invertible $\sO_X$-modules; its group of isomorphism classes is denoted by $\Pic(X)$. The category of locally free $\sO_X$-modules is denoted by $\mathbf{Vect}_X$, and the category of reflexive $\sO_X$-modules by $\sO_X$-$\mathbf{Rflx}$. 

By a semilocal scheme we mean an affine scheme with semilocal coordinate ring.

\epp

\subsection*{Acknowledgements}
We would like to express our gratitude to K\k{e}stutis \v{C}esnavi\v{c}ius for his constant encouragement throughout this project.
We are also indebted to Matthew Morrow and Colliot-Thélène for proposing the Grothendieck--Serre conjecture on smooth schemes over semilocal Pr\"ufer rings during the defence of the first author.
We are grateful to K\k{e}stutis \v{C}esnavi\v{c}ius for his help in removing the assumptions on finite residue fields in our original formulation of the main \Cref{main thm}. Additionally, after an earlier version of this paper was completed, K\k{e}stutis \v{C}esnavi\v{c}ius kindly sent us his note, which contained a sketch of a different proof of the Noetherian case of \Cref{main thm}.
On several occasions over the past few months, we had conversations with K\k{e}stutis \v{C}esnavi\v{c}ius, Arnab Kundu, and Shang Li about various aspects of this article, for which we are grateful.
Finally, we acknowledge Jiandi Zou for his useful suggestions.
 This project has received funding from the European Research Council (ERC) under the European Union's Horizon 2020 research and innovation programme (grant agreement No. 851146), the grant 075-15-2022-289, and the excellent environment for research provided by the Harbin Institute of Technology and the Euler International Mathematical Institute.


\section{Purity for torsors under reductive groups}
\label{sect-purity of reductive torsors}
 
 The main result of this section is \Cref{extend generically trivial torsors}, which extends the domain of definition of a generically trivial torsor under a reductive group scheme on a scheme smooth over a Pr\"ufer base. It will be used in the proof of \Cref{torsors-Sm proj base}. We first recall the geometric facts about schemes over Pr\"ufer bases that are needed below. Much of the material on coherence and reflexivity is standard and appears in Gabber--Ramero's work \cite{GR18}; we recall the relevant statements here.

Recall that a ring $A$ is \emph{coherent} if every finitely generated ideal of $A$ is finitely presented. A scheme is \emph{locally coherent} if it admits an affine open cover by spectra of coherent rings. 
On a locally coherent scheme $X$, coherent $\sO_X$-modules are exactly finitely presented $\sO_X$-modules. See, \emph{e.g.}, \cite{Gla89}*{Corollary~2.2.2 and Theorem~2.3.2}.

\blemt \label{geom}
 Let $X$ be a scheme that is flat and locally of finite type over an integral Pr\"ufer scheme $S$.
\benumr
  \item \label{coherence of X}The scheme $X$ is locally of finite presentation and locally coherent.
  \item \label{coherence of O}For every point $x\in X$, the local ring $\sO_{X,x}$ is coherent.
  \item\label{geo-i} If $X$ is irreducible, then all the nonempty $S$-fibres have the same dimension;
  \item\label{geo-iii} If $s\in S$ and $\sO_{X_s,\xi}$ is reduced for a maximal point $\xi\in X_s$, then the local ring $\sO_{X,\xi}$ is a valuation ring and the extension $\sO_{S,s}\hra \sO_{X,\xi}$ induces an isomorphism of value groups.
  \item\label{limitarg} If $S=\Spec R$ is semilocal affine with fraction field $K$ and $X$ is of finite type over $S$, then, for some finitely generated subfield $K_0\subset K$, with $R_0\ce R\cap K_0$ and $S_0\ce \Spec R_0$, the ring $R_0$ is a semilocal Pr\"ufer domain of finite Krull dimension, the inclusion $R_0\subset R$ induces $S\to S_0$, and $X\simeq X_0\times_{S_0}S$ for a flat finite type $S_0$-scheme $X_0$.
 \eenum
\elemt
\bpf
For \ref{coherence of X}--\ref{coherence of O}, see \cite{GL24}*{Lemma~3.2.1}. For \ref{geo-i}, see \cite{EGAIV3}*{lemme~14.3.10}.
For \ref{geo-iii}, see \cite{MB22}*{théorème~A}. 
For \ref{limitarg}, write $K=\Frac R$. Since $R$ is a semilocal Pr\"ufer domain, $R=\bigcap_{\fm\in\operatorname{Max}R}R_{\fm}$ inside $K$. For a finitely generated subfield $K_0\subset K$, put
\[
   \tst R_0\ce R\cap K_0=\bigcap_{\fm\in\operatorname{Max}R}(R_{\fm}\cap K_0).
\]
The criterion in \cite{Bou98}*{chapitre~VI, \S1, \textnumero~4, exemple~(2)} shows that each $R_{\fm}\cap K_0$ is a valuation ring of $K_0$. Since $K_0$ is finitely generated over its prime field, Abhyankar's inequality \cite{Abh56v}*{Theorem~1, p.~330} and \cite{Bou98}*{chapitre~VI, \S10, \textnumero~2, proposition~3 et corollaire} imply that these valuation rings have finite rank. Thus, by \cite{Bou98}*{chapitre~VI, \S7, \textnumero~1, propositions~1 et~2}, the rings $R_0$ form a filtered system of semilocal Pr\"ufer domains of finite Krull dimension with colimit $R$. By \ref{coherence of X}, the finite type $R$-scheme $X$ is finitely presented. 
The approximation \cite{EGAIV3}*{th\'eor\`eme~8.8.2} then produces, after enlarging $K_0$, a finite type $R_0$-scheme $X_0$ with $X\simeq X_0\times_{S_0}S$, and flatness descends by \cite{EGAIV3}*{proposition~11.2.6}.
\epf

\bpp[Setup]\label{setup}
For convenience, we fix the following notations.
Let
\begin{itemize}
\item $R$ be a semilocal Pr\"ufer domain with spectrum $S$ and generic point $\eta\in S$;
\item $f\colon X\ra S$ be a smooth, finite type morphism of schemes;
\item $Z\subset X$ be a closed subset such that $j\colon X\backslash Z\hra X$ is quasi-compact with the following property
\[
\q \text{$\codim(Z_{\eta},X_{\eta})\geq 2$\q and \q $\codim(Z_s, X_s)\geq 1$ for $\eta \neq s\in S$;}
\]
\item $G$ be a reductive $X$-group scheme.
\end{itemize}
\epp
Note that when the relative dimension of $f$ is at most $1$, the condition on $Z$ above implies that $Z_{\eta}=\emptyset$.
On the locally coherent schemes considered here, a coherent $\sO_X$-module $\sF$ is reflexive if and only if it is \emph{finitely copresented}, that is, there are positive integers $m$ and $n$ such that $\sF$ fits into the following exact sequence
\[
   0\ra \sF\ra \sO_X^{\oplus m}\xrightarrow{h} \sO_X^{\oplus n}.
\]
We will use the following consequence of \cite{GR18}*{Lemma~11.3.3}. For reduced locally coherent schemes, the dual of a finitely presented module is reflexive.
The rank condition in that lemma is local on the scheme and, for finitely presented modules, may be checked at maximal points.
Indeed, for a reflexive $\sO_X$-module $\sF$, its dual $\sF^{\vee}$ is finitely presented. 
Locally, after choosing a finite presentation of $\sF$, the dual $\sF^{\vee}$ is the kernel of a map between finite free modules, hence is coherent by \SP{01BY}, and therefore finitely presented. Thus such a copresentation can be obtained as the dual of a finite presentation of $\sF^{\vee}$; conversely, if $\sF$ admits such a copresentation, then 
\[
   \x{$\sF=\ker(h)=(\coker(h^{\vee}))^{\vee}$\q is $\sO_X$-reflexive.}
\]

We then present Gabber--Ramero's extension  for reflexive modules \cite{GR18}*{Proposition~11.3.8}.

\bpropt\label{GR18}
The restrictions induce the following equivalences, with quasi-inverses given by $j_*$:
\begin{equation}\label{restriction}
\sO_{X}\x{-}\mathbf{Rflx}\isoto \sO_{X\backslash Z}\x{-}\mathbf{Rflx}
\q and \q 
\mathbf{Pic}_X\isoto \mathbf{Pic}_{X\backslash Z}.
\end{equation}
\epropt
\bpf 
The problem is Zariski-local on $X$, so we may assume that $X$ is affine.
Note that $X\backslash Z$ is quasi-compact.
First, we show that the restrictions considered are essentially surjective.
Given a reflexive $\sO_{X\backslash Z}$-module $\sF$, its dual $\sF^{\vee}$ is a coherent $\sO_{X\backslash Z}$-module, so, by \cite{GR18}*{Lemma~10.3.24(ii)}, $\sF^{\vee}$ extends to a finitely presented $\sO_X$-module $\wt{\sF}$. 
By the consequence of \cite{GR18}*{Lemma~11.3.3} above, the dual $(\wt{\sF})^{\vee}$ is reflexive. Since duals of finitely presented modules commute with restriction to the open $X\backslash Z$, this reflexive module restricts to $(\sF^\vee)^\vee\simeq \sF$. 
Next, to deduce full faithfulness, we let $\sF$ and $\sG$ be reflexive $\sO_X$-modules and show that restriction on $U\ce X\backslash Z$ induces 
\[
   \Hom_{\sO_X}(\sF,\sG)\isoto \Hom_{\sO_U}(\sF|_{U},\sG|_{U}).
\]
Since $\sF$ is finitely presented, we are reduced to the case when $\sF$ is free.
Since $\sG$ is reflexive and hence finitely copresented, it suffices to assume that $\sG$ is free.
The problem is reduced to showing that 
\[
   \x{$\sO_X(X)\isoto \sO_X(U)$\q is an isomorphism.}
\]
Consider the functor $\sH^0_Z\colon \sF\mapsto \ker(\sF\to j_{\ast}j^{\ast}\sF)$ on $\sO_X$-modules.  
As $\sH^0_Z$ is left exact, it has right derived functors, denoted by $\sH^i_Z$. 
We have an isomorphism $R^0\sH^0_Z\cong \sH^0_Z$ and the  exact sequence 
\[
   0\ra \sH^0_Z(\sO_X)\ra \sO_X\ra j_{\ast}j^{\ast}\sO_X \ra \sH^1_Z(\sO_X).
\]
The problem is local; by induction we may assume that $Z=\{z\}$ is the closed point of $X=\Spec \sO_{X,z}$. 
The vanishing $\sH^0_Z(\sO_X)=0$ is immediate. If $z$ lies over a non-generic point of $S$, a section supported at $z$ dies on the generic fibre, hence is zero because $\sO_X$ is $R$-torsion-free by flatness. If $z$ lies over $\eta$, then the local ring is $\sO_{X_{\eta},z}$, which is regular, so it has no nonzero element supported at its closed point.
Thus, it remains to show that $\sH^1_Z(\sO_X)=0$.
It remains to produce an $\fm_z$-regular sequence of length at least two.
If $z\in X_\eta$, then the bound $\codim(Z_\eta,X_\eta)\geq 2$ and the regularity of the smooth generic fibre produce an $\fm_z$-regular sequence in $\sO_X=\sO_{X_\eta,z}$ of length $2$.
If $z$ lies over a non-generic $s\in S$, put $A\ce\sO_{S,s}$ and $B\ce\sO_{X,z}$.
The inequality $\codim(Z_s,X_s)\geq 1$ implies $\dim(B\otimes_A k_s)>0$.
Since $B\otimes_A k_s\simeq \sO_{X_s,z}$ is regular local, choose a nonzerodivisor $\bar f$ in its maximal ideal and lift $\bar f$ to an element $f\in \fm_B$.
By \cite{EGAIV3}*{proposition~11.3.7}, the element $f$ is $B$-regular and $B/fB$ is $A$-flat.
Choose a nonzero $a\in\fm_A$, so multiplication by $a$ is injective on $B/fB$.
Thus $(f,a)$ is an $\fm_z$-regular sequence in $B$.
In either case, choose an $\fm_z$-regular sequence $(f_1,\ldots,f_r)$ with $r\geq 2$.
A regular element $f_1\in \fm_z$ induces a short exact sequence 
\[
   0\ra \sO_X\overset{\times f_1}{\ra}\sO_X\ra \sO_X/f_1\ra 0,
\]
which leads to a long exact sequence of $\sH^i_Z(\sO_X)$.
Note that the $\sH^i_Z(-)$ are supported at $z$, so, iterating along the regular sequence $(f_1,\ldots,f_r)$, we deduce that $\sH^i_Z(\sO_X)=0$ when $i<r$.
In particular, we have $\sH^0_Z(\sO_X)=\sH^1_Z(\sO_X)=0$, so the restriction functor on reflexive sheaves is fully faithful. This finishes the proof for the reflexive modules.

It remains to show that, for an invertible $\sO_{X\backslash Z}$-module  $\sL$, the unique reflexive extension $j_{\ast}\sL$ is $\sO_X$-invertible.
Recall \cite{KM76}*{Theorem~1} that there exists the determinant functor 
\[
   \det \colon  D(\sO_{X}\x{-}\mathbf{Mod})^{\ast}_{\x{perf}}\ra \mathbf{gr.Pic}(X)
\]
from the groupoid of perfect complexes of $\sO_X$-modules to the category of graded invertible $\sO_X$-modules. The functor $\det$ commutes with arbitrary base change.
The reflexive module $j_{\ast}\sL$ is finitely presented.
At a point $x\in X$ with image $s\in S$, apply \cite{GR18}*{Proposition~11.4.1(iii)} to $X_{\sO_{S,s}}\to \Spec\sO_{S,s}$ and to the stalk $(j_{\ast}\sL)_x$.
This shows that the reflexive module has finite projective dimension.
Since it is finitely presented and has finite projective dimension, $j_{\ast}\sL$ is perfect locally on $X$.
Thus the invertible $\sO_X$-module $\det(j_{\ast}\sL[0])$ is well-defined. We have
\[
   \begin{aligned}
   \det(j_{\ast}\sL[0])
   &\isoto j_{\ast}j^{\ast}\det(j_{\ast}\sL[0])
   \isoto j_{\ast}\det(j^{\ast}j_{\ast}\sL[0])  \\
   &\isoto j_{\ast}\det(\sL[0])
   \isoto j_{\ast}\sL.
   \end{aligned}
\]
Here, the first isomorphism follows from the assertion for reflexive modules, and the second isomorphism uses the base-change property of $\det$ together with $j^{\ast}j_{\ast}\sL\simeq \sL$.
This shows that $j_{\ast}\sL\simeq \det(j_{\ast}\sL[0])$ is $\sO_X$-invertible, and hence proves the Picard equivalence.
\epf

The following result provides an analogue of the purity theorem of Colliot-Th\'el\`ene and Sansuc for reductive group torsors over two-dimensional regular Noetherian schemes (\cite{CTS79}*{th\'eor\`eme 6.13}).
\bthmt\label{purity for rel. dim 1}
If $f$ has relative dimension $1$, then restriction induces the following equivalence 
\begin{equation}\label{restriction of G-torsors}
  \mathbf{Tors}(X_{\et},G) \isoto \mathbf{Tors}((X\backslash Z)_{\et},G)
\end{equation}
of categories of $G$-torsors.
Moreover, restriction induces an equivalence $\mathbf{Vect}_X\isoto \mathbf{Vect}_{X\backslash Z}$.
In particular, passing to isomorphism classes of objects, we have a bijection $$H^1_{\et}(X,G)\simeq H^1_{\et}(X\backslash Z,G).$$
\ethmt
\bpf 
We begin with the case when $G=\GL_n$. 
\Cref{GR18} applied  to the reflexive module $\sO_X$ yields $\sO_X\simeq j_{\ast}\sO_{X\backslash Z}$. 
It follows that $Y(X)=Y(X\backslash Z)$ for every $X$-affine scheme $Y$.
In particular, this applies to the $X$-affine scheme  $Y\ce {\Isom}_X(\sE_1,\sE_2)$, where $\sE_1$ and $\sE_2$ are vector bundles over $X$.
This proves the full faithfulness of the restriction $\mathbf{Vect}_X\ra \mathbf{Vect}_{X\backslash Z}$. 
For the essential surjectivity, by \Cref{GR18}, every vector bundle $\sE$ on $X\backslash Z$ extends to a reflexive $\sO_X$-module $\wt{\sE}\ce j_{\ast}\sE$. 
We show that $\wt{\sE}$ is locally free.
As $\wt{\sE}$ is finitely copresented, it is $R$-torsion-free, hence $R$-flat.
For any $x\in X$ with image $s\in S$, the fibre $(\wt{\sE})_s$ is a finite torsion-free module on the regular curve $X_s$, hence is locally free.
Combining \cite{EGAIV3}*{proposition~11.3.7} with Nakayama's lemma lifts a basis of $(\wt{\sE})_s$ at $x$ to a basis of $\wt{\sE}$ around $x$.
Thus $\wt{\sE}$ is a vector bundle.

For a general $G$, the full faithfulness of $\mathbf{Tors}(X_{\et},G) \isoto \mathbf{Tors}((X\backslash Z)_{\et},G)$ can be obtained similarly, since the isomorphism scheme for two $G$-torsors is again $X$-affine by descent.
For the essential surjectivity, we start by picking a $G$-torsor $\cP$ on $X\backslash Z$ and seek to show that $\cP$ extends to a $G$-torsor on $X$. 
Note that the codimension condition on $Z$ is preserved under \'etale base change. By gluing, we are reduced to showing that, \emph{\'etale locally} on $X$, $\cP$ extends over $X$.
We may assume that $X$ is affine and that $G\subset \GL_{n,X}$. Consider the following commutative diagram with exact rows
\[
\begin{tikzcd}
{(\mathrm{GL}_{n,X}/G)(X)} \arrow[r] \isoarrow{d} & {H^1_{\et}(X,G)} \arrow[r] \arrow[d] & {H^1_{\et}(X,\mathrm{GL}_{n,X})} \arrow[d] \\
{(\mathrm{GL}_{n,X}/G)(X\backslash Z)} \arrow[r]           & {H^1_{\et}(X\backslash Z,G)} \arrow[r]           & {H^1_{\et}(X\backslash Z,\mathrm{GL}_{n,X})},
\end{tikzcd}
\]
where the bijectivity of the left vertical arrow follows since $\GL_{n,X}/G$ is $X$-affine (see \cite{Alp14}*{Theorem~9.4.1}) and sections of $X$-affine schemes extend uniquely across $Z$.
By the vector bundle case above, we may replace $X$ by an affine open cover to assume that the induced $\GL_{n,X\backslash Z}$-torsor $\cP\times^{G_{X\backslash Z}}\GL_{n, X\backslash Z}$ is trivial. 
A diagram chase implies that there exists a $G$-torsor $\cQ$ on $X$ such that $\cQ|_{X\backslash Z}\simeq \cP$, as desired. 
\epf 
The following local variant of \Cref{purity for rel. dim 1} is a non-Noetherian counterpart of Colliot-Thélène--Sansuc's result \cite{CTS79}*{théorème~6.13}.

\bthmt \label{extends across codim-2 points}
Assume that $R$ is a finite rank valuation ring.
For a point $x\in X$ that is
\benumr
\item\label{g-tor-local-gen} either $x\in X_{\eta}$ with $\dim \sO_{X_{\eta},x} =2$, or
\item\label{g-tor-local-ngen} $x\in X_s$ for some non-generic $s\in S$ with $\dim \sO_{X_s,x} =1$,
\eenum
put $X_x\ce \Spec \sO_{X,x}$, $U_x\ce X_x\backslash\{x\}$, and $G_x\ce G_{X_x}$.
Then the restriction functor
\[
   \mathbf{Tors}((X_x)_{\et},G_x)\ra \mathbf{Tors}((U_x)_{\et},G_{U_x})
\]
is an equivalence.
Moreover, restriction induces an equivalence
\[
   \mathbf{Vect}_{X_x}\isoto \mathbf{Vect}_{U_x}.
\]
In particular, every $G_{U_x}$-torsor on $U_x$ extends uniquely to a $G_x$-torsor on $X_x$.
\ethmt
\bpf 
Let $X_0\ce X_x$, let $U\ce U_x$, let $j\colon U\hra X_0$ be the open immersion, and let $G_0\ce G_x$.
The open $U$ is quasi-compact. In case \ref{g-tor-local-gen}, this is the punctured spectrum of the two-dimensional Noetherian regular local ring $\sO_{X_{\eta},x}$. In case \ref{g-tor-local-ngen}, the finite-rank assumption on $R$ makes $X$ topologically Noetherian, hence so is $X_0$.

The part of \Cref{GR18} needed below says that sections of $X_0$-affine schemes extend uniquely from $U$ to $X_0$.
The local cohomology argument of that proposition, applied to the pair $(X_0,\{x\})$, proves
\[
   H^0_{\{x\}}(X_0,\sO_{X_0})=H^1_{\{x\}}(X_0,\sO_{X_0})=0.
\]
In case \ref{g-tor-local-gen}, the local ring $\sO_{X_0,x}$ is regular of dimension $2$.
In case \ref{g-tor-local-ngen}, choose a regular element $\bar f$ in the maximal ideal of the regular one-dimensional local ring $\sO_{X_s,x}$ and lift it to $f\in \sO_{X_0,x}$.
By \cite{EGAIV3}*{proposition~11.3.7}, the element $f$ is regular and $\sO_{X_0,x}/f$ is $\sO_{S,s}$-flat.
For any nonzero $a\in\mathfrak m_s$, the pair $(f,a)$ is then a regular sequence.
Thus in both cases the depth at $x$ is at least $2$, which proves the displayed vanishing.
Consequently, $\sO_{X_0}\simeq j_{\ast}\sO_U$.
It follows that for every $X_0$-affine scheme $\mathcal W=\Spec_{X_0}\mathcal B$,
\[
   \mathcal W(X_0)=\Hom_{\sO_{X_0}\text{-alg}}(\mathcal B,\sO_{X_0})
   \isoto
   \Hom_{\sO_U\text{-alg}}(\mathcal B|_U,\sO_U)=\mathcal W(U).
\]

This proves full faithfulness for both torsors and vector bundles. Indeed, the corresponding isomorphism schemes are $X_0$-affine, so their sections over $X_0$ and over $U$ agree.

We next prove the essential surjectivity for vector bundles.
For a vector bundle $\sE_U$ on $U$, the reflexive extension in \Cref{GR18} produces a finite reflexive $\sO_{X_0}$-module $\wt{\sE}$ with $\wt{\sE}|_U\simeq \sE_U$.
In case \ref{g-tor-local-gen}, the ring $\sO_{X_0,x}$ is two-dimensional regular, so the finite reflexive module $\wt{\sE}$ is free at $x$ by the Auslander--Buchsbaum formula.
In case \ref{g-tor-local-ngen}, the criterion \cite{GR18}*{Proposition~11.4.1(iii)} shows that the projective dimension at $x$ is $0$, because $\dim\sO_{X_s,x}=1$.
Thus $\wt{\sE}$ is locally free.
Together with full faithfulness, this proves $\mathbf{Vect}_{X_0}\simeq \mathbf{Vect}_{U}$.

It remains to prove the essential surjectivity for torsors.
In case \ref{g-tor-local-gen}, this is Gabber's purity theorem for reductive torsors over the punctured spectrum of a two-dimensional regular local ring, see \cite{Gab81}*{\S1, Lemma~1}.
Assume from now on that \ref{g-tor-local-ngen} holds.

Now let $\cP$ be a $G_U$-torsor on $U$.
Since extension can be checked \'etale locally on $X_0$ and uniqueness has already been proved, we may work after an \'etale localization at $x$ and assume that $G_0$ is a closed subgroup of some $\GL_{n,X_0}$.
Let $\cP_{\GL}\ce\cP\times^{G_U}\GL_{n,U}$, viewed as a torsor under $\GL_{n,U}$.
The torsor $\cP_{\GL}$ corresponds to a vector bundle on $U$.
By the equivalence $\mathbf{Vect}_{X_0}\simeq \mathbf{Vect}_{U}$ proved above, this vector bundle extends to $X_0$.
Since $X_0$ is local and affine, the extension is free, so $\cP_{\GL}$ is trivial.
Choosing such a trivialization identifies $\cP$ with a section over $U$ of the $X_0$-affine scheme $\GL_{n,X_0}/G_0$.
The section extends uniquely to $X_0$ by the affine extension property above.
Pulling back $\GL_{n,X_0}\to \GL_{n,X_0}/G_0$ along the extended section produces a $G_0$-torsor on $X_0$ whose restriction to $U$ is $\cP$.
The local extensions glue by full faithfulness.
\epf

The following corollaries of the purity theorem above allow us to extend torsors under reductive group schemes outside a closed subset of higher codimension, which is crucial for our proof of the Grothendieck--Serre conjecture for constant reductive groups in \Cref{G-S for constant reductive gps}.
\bpropt \label{extension torsors}
For a  $G$-torsor $\cP$ on $X\backslash Z$, there are
\benumr 
\item a closed subset $Z\pr\subset Z$ satisfying $\codim(Z\pr_{\eta}, X_{\eta})\geq 3$ and  $\codim(Z\pr_s,X_s)\geq 2$\,\, for all $s\in S$, and 
\item a $G$-torsor $\cQ$ on $X\backslash Z\pr$
\eenum
such that $\cP\simeq \cQ|_{X\backslash Z}$.
\epropt
\bpf 

A smooth morphism is locally quasi-separated.
Since $f$ is finite type, it is quasi-compact.
Thus $f$ is quasi-separated, so it is of finite presentation.
We may now apply the limit argument \Cref{geom}\ref{limitarg} and reduce to the case when $R$ has finite Krull dimension, so $X$ is topologically Noetherian. 
In particular, there are at most finitely many points $x\in Z$ satisfying either assumption \ref{g-tor-local-gen} or \ref{g-tor-local-ngen} of \Cref{extends across codim-2 points}. 
Thus, by gluing and spreading out finitely many times, $\cP$ extends across these points.
Eventually, we find a desired closed subset $Z\pr\subset X$ such that $\cP_{X\backslash Z}$ extends to a $G$-torsor $\cQ$ over $X\backslash Z\pr$.
\epf

\bcort \label{extend generically trivial torsors}
Let $\mathbf{x}\subset X$ be a finite subset contained in an affine open. 
Let $\cO$ be the semilocalization of $X$ at $\mathbf{x}$. 
Pick a nonzero $r\in \cO$.
Then, every generically trivial $G$-torsor on $\cO[\f{1}{r}]$ extends to a $G$-torsor over $X\backslash Y$ for a closed subset $Y\subset X$ satisfying the following conditions
\[
\text{$\codim(Y_{\eta}, X_{\eta})\geq 3$ \q and\q  $\codim(Y_s,X_s)\geq 2$\,\, for all $\eta\ne s\in S$.}
\]
\ecort
\bpf 
Let $\cP$ be a generically trivial $G$-torsor as in the statement.
Choose an affine open neighbourhood $V=\Spec B\subset X$ of $\mathbf{x}$.
At the start of the proof of \Cref{extension torsors} we checked that $f$ is quasi-separated and of finite presentation.
As $S$ is affine, $X$ is quasi-separated.
Thus $V\hra X$ is quasi-compact.
Being an open immersion, it is of finite presentation, so $V$ is of finite presentation over $S$.
After replacing $V$ by a principal open containing $\mathbf{x}$, the element $r\in\cO$ is represented by some $r_B\in B$, and the finitely presented torsor $\cP$ descends to a $G$-torsor $\cP_U$ on $U\ce D(r_B)\subset V$.
The data $X,V,U,G,\cP_U$, together with the open immersions $U\subset V\subset X$, are then of finite presentation over $S$.
By \Cref{geom}\ref{limitarg}, after enlarging the finite-dimensional $S_0$, these data descend to $S_0$.
Any extension constructed over $S_0$ pulls back to an extension over $S$.
Hence we may assume that $R$ has finite Krull dimension, so $X$ is topologically Noetherian.
It remains to enlarge $U$ and extend $\cP_U$ to ensure that $Z\ce X\backslash U$ satisfies the assumptions of \Cref{extension torsors}.
Assume that $Z$ contains a point $z$ such that one of the following holds
\benumr
\item $z\in X_{\eta}$ and $\dim \sO_{X,z}=1$ (so $U\cap \Spec \sO_{X,z}$ is a maximal point of $X$);
\item $z$ is a maximal point of $X_s$ for a non-generic point $s\in S$ (so $U\cap \Spec \sO_{X,z}$ is the spectrum of a valuation ring).
\eenum 
By the Grothendieck--Serre theorem over valuation rings \cite{Guo24}, $\cP_U|_{\Spec \sO_{X,z}\cap U}$ is trivial.
Thus, we can glue $\cP_U$ with the trivial $G$-torsor on a small neighbourhood of $z\in Z$ to extend $\cP_U$ across such $z$. 
As the topological of $X$ is Noetherian, $Z$ has finitely many points of types (i) and (ii). 
Therefore, iteratively extending $U$ and $\cP_U$, we may assume that $Z$ has no points of type (i) or (ii), whence \Cref{extension torsors} applies.
\epf

\section{Low dimensional cohomology of groups of multiplicative type}
\label{low dim cohomology of mult gp}
The primary result of this section is the Grothendieck--Serre theorem for groups of multiplicative type over smooth schemes over Pr\"ufer bases (\Cref{G-S type results for mult type}), with a particular emphasis on tori. 
This theorem serves as a foundational aspect of the Grothendieck--Serre conjecture over semilocal Pr\"ufer rings (\Cref{G-S over semi-local prufer}).
To establish \Cref{G-S type results for mult type}, we start with the purity \Cref{purity for gp of mult type} and then apply a criterion for the descent of line bundles (\Cref{non-Noeth-pullback line bundle}), employing an approach distinct from the classical one used in \cite{CTS87}*{4.1--4.3}.
\bthm \label{purity for gp of mult type}
Under the Setup \ref{setup}, let $M$ be a finite type $X$-group of multiplicative type and $U\ce X\backslash Z$.
\benumr
  \item \label{mult-paraf} We have
\[
   H^i_{\fppf}(X,M)\simeq H^i_{\fppf}(U,M)\q \,\text{for } i=0,1 \q \text{and} \q\, H^2_{\fppf}(X,M)\hra H^2_{\fppf}(U,M).
\]
    \item  \label{gl-H1-H2-flasque} Assume that $M=T$ is an $X$-torus and  $Y\subset X$ is a nowhere dense finitely presented  closed subscheme.
    If $T_{\sO_{X,y}}$ is a flasque torus\footnote{For a torus over a connected base $W$, choose a connected finite \'etale Galois cover $\wt{W}\to W$ splitting the torus, with Galois group $\Gamma$, and let $\Lambda_T$ be the resulting cocharacter $\Gamma$-lattice. Following \cite{Ces22}*{Appendix~A.2}, the torus is flasque if $\Ext^1_{\bZ[\Gamma]}(P,\Lambda_T)=0$ for every quasi-trivial $\Gamma$-lattice $P$, that is, for every $\Gamma$-lattice with a $\Gamma$-stable $\bZ$-basis. This condition is independent of the chosen splitting cover and is checked componentwise on the base.} for every $y\in Y$ for which the local ring $\sO_{X,y}$ is a valuation ring, then we have
  \[
\text{$H^1_{\et}(X,T) \twoheadrightarrow H^1_{\et}(X\backslash Y,T)$ \q and \q $H^2_{\et}(X,T)\hra H^2_{\et}(X\backslash Y,T)$.}
\]
In particular, if $K(X)$ denotes the total ring of fractions of $X$, then we have
 \[
\text{$\Pic(X)\twoheadrightarrow \Pic(X\backslash Y)$ \q and \q $\mathrm{Br}(X)\hra \mathrm{Br}(K(X))$.}
\]
\eenum
\ethm
\bpf
By \Cref{geom}\ref{limitarg} and limit formalism for the fixed cohomology groups that occur in the two assertions, the relevant restriction maps are obtained by base change over a semilocal Pr\"ufer domain of finite Krull dimension. Thus, before passing to local cohomology sheaves, we may assume that $R$ has finite Krull dimension. In particular, the topological space of $X$ is Noetherian.

For the assertion \ref{mult-paraf}, by the local cohomology exact sequence for the pair $(X,Z)$ and the sheaf $M$,
the statement is equivalent to the vanishing $H^i_Z(X,M)=0$ for $0\le i \le 2$. By the spectral sequence in \cite{CS24}*{Lemma 7.1.1}, it suffices to show the vanishing of $\cH^q_Z(M)$, the sheafification of the presheaf $X\pr\mapsto H^q_{Z\pr}(X\pr, M)$ on $X_{\et}$, where $Z\pr\ce Z\times_X X\pr$. 
This problem is \'etale-local on $X$, so we may assume that $X$ is affine, and $M$ splits as a product of $\bG_m$ or $\mu_n$.
As $\mu_n=\ker(\bG_m\overset{\times n}{\ra}\bG_m)$, it remains to consider the case when $M=\bG_m$.
Recall the coniveau spectral sequence \cite{Gro68c}*{\S10.1}
 \[
 \text{{$E^{pq}_1$}$=\bigoplus_{z\in Z^{(p)}} H^{p+q}_{\{z\}}(\bG_m)\Rightarrow H^{p+q}_Z(X,\bG_m)$, \q where \q $H^{p+q}_{\{z\}}(\bG_m)\ce \operatorname*{colim}_{W\ni z}\,H^{p+q}_{\overline{\{z\}}\cap W}{(W,\bG_m)}$}
 \]
 and $W$ runs over the open neighbourhoods of $z$ in $X$. As the topological space of $X$ is Noetherian, we can identify $ H^{p+q}_{\{z\}}(\bG_m)$ with $H^{p+q}_{\{z\}}(\sO_{X,z},\bG_m)$.
   Therefore, it is enough to show $H^i_{\{z\}}(\sO_{X,z},\bG_m)=0$ for $z\in Z$ and $0\le i \le 2$, which will be settled by \Cref{local cohomology of tori}\ref{paraf}.

For the assertion \ref{gl-H1-H2-flasque}, the local cohomology exact sequence reduces us to showing that $H^2_Y(X,T)=0$. 
By the coniveau spectral sequence, it is enough to prove that $H^2_{\{y\}}(\sO_{X,y},T)=0$ for every $y\in Y$.
Since $Y$ contains no irreducible component of $X$, if $y\in X_{\eta}$, then $\dim\sO_{X_{\eta},y}\neq 0 $.
If $y\in X_{\eta}$ with $\dim\sO_{X_{\eta},y}\ge 2$, or if $y\in X_s$ for some non-generic $s\in S$ with $\dim\sO_{X_s,y}\ge 1$, then $\overline{\{y\}}$ satisfies the codimension conditions in Setup~\ref{setup}.
Thus \Cref{local cohomology of tori}\ref{paraf} proves the desired vanishing.
It remains to treat the codimension-$1$ points of $X_{\eta}$ and the maximal points of the non-generic fibres.
In the first case $\sO_{X,y}$ is a discrete valuation ring. 
In the second case, \Cref{geom}\ref{geo-iii} shows that $\sO_{X,y}$ is a valuation ring with the same value group as $\sO_{S,s}$.
All these valuation rings are of finite rank.
Therefore, \Cref{local cohomology of tori}\ref{val} applies by the flasqueness assumption on $T_{\sO_{X,y}}$.
\epf

\blem\label{local cohomology of tori}
Under the Setup \ref{setup}, assume that $R$ is a finite-rank valuation ring, and let $x\in X$ and $T$ be an $\sO_{X,x}$-torus.
\benumr
\item \label{paraf} If $x\in Z$, then we have
    \[
    \text{$H^i_{\{x\}}(\sO_{X,x},T)=0$ \q for \, $0\le i \le 2$. }
    \]
\item \label{val} If $\sO_{X,x}$ is a finite-rank valuation ring that is not a field and $T$ is flasque, then we have
     \[
    H^2_{\{x\}}(\sO_{X,x},T)=0.
    \]
\eenum
\elem
\brem
Without assuming that $x\in Z$, the assertion \ref{paraf} still holds for any $x\in X$ such that
\begin{equation}\label{pt-cond}
   \x{either $x\in X_{\eta}$ with $\dim \sO_{X_{\eta},x} \ge 2$, \qq or $x\in X_s$ with $\dim \sO_{X_s,x} \ge 1$.}
\end{equation}
To see this, by \Cref{geom}\ref{geo-i}, the closure $\overline{\{x\}}$ is $R$-fibrewise of codimension $\ge 2$ (resp., $\ge 1$) in $X_{\eta}$ (resp., in $X_s$), so it is a valid choice for $Z$.
\erem
\begin{proof}[Proof of \Cref{local cohomology of tori}]
For the assertion \ref{val}, put $V\ce\sO_{X,x}$.
Since $V$ is a finite-rank valuation ring that is not a field, \cite{Guo24}*{Lemma~2.6} proves that $H^2_{\{x\}}(V,T)=0$.

For the assertion \ref{paraf}, it suffices to check the stalks after strict Henselization. 
Fix a geometric point $\bar{x}$ over $x$, put $A\ce \sO_{X,\bar{x}}^{\sh}$, and set $P_A\ce \Spec A\backslash\{\bar{x}\}$. The torus $T_A$ is split, so it is enough to treat $T=\bG_m$. Consider the cofiltered category of affine pointed \'etale neighbourhoods $(\cW,\bar{x}\to y\in \cW)$ of $\bar{x}$, compatible with the chosen geometric point, and put $B_{\cW}\ce \sO_{\cW,y}$. For every such $\cW$, applying \Cref{GR18} to a small affine neighbourhood of $y$ and then localizing at $y$ yields
\[
H^i_{\et}(\Spec B_{\cW},\bG_m)\simeq H^i_{\et}(\Spec B_{\cW}\backslash\{y\},\bG_m) \q \text{ for } \q i=0,1.
\]
Here $\Spec B_{\cW}\backslash\{y\}$ is quasi-compact. Since the affine map $\Spec A\to \Spec B_{\cW}$ is quasi-compact, the punctured spectrum $P_A$ is also quasi-compact, and the above isomorphisms pass to the filtered limit and yield
\[
H^i_{\et}(\Spec A,\bG_m)\simeq H^i_{\et}(P_A,\bG_m) \q \text{ for } \q i=0,1.
\]
Examining the local cohomology exact sequence for the pair $(\Spec A,\{\ov{x}\})$ and the split torus $T_A\simeq \mathbb{G}_{m}^{\dim T}$, we deduce that
\[
   H_{\{\ov{x}\}}^i(\sO_{X,\bar{x}}^{\sh},T)=0 \q \text{ for } \q 0\le i\le 2.
\]
By the local-to-global spectral sequence \cite{SGA4II}*{exposé~V, proposition~6.4}, we complete the proof.
\end{proof}

The following lemma for descending line bundles that are trivial on generic fibres is useful for \Cref{G-S type results for mult type}.
\blem \label{non-Noeth-pullback line bundle}
Let $\phi\colon X\to Y$ be a morphism of schemes. Let $\sL$ be an invertible $\sO_X$-module. If
\begin{itemize}
  \item [(1)] $Y$ is quasi-compact, quasi-separated, integral, and normal,
  \item [(2)] there exist a smooth projective morphism $\overline{\phi}\colon \overline{X}\to Y$, with geometrically integral fibres, and a quasi-compact open immersion $X\hookrightarrow \overline{X}$ over $Y$, and
  \item [(3)] $\sL$ is trivial when restricted to the generic fibre of $\phi$,
\end{itemize}
then $\sL\simeq \phi^*\sN$ for some invertible $\sO_Y$-module $\sN$.
\elem
\bpf
When $Y$ is Noetherian, the lemma follows from a much more general result \SP{0BD6}, where (2) can be replaced by the weaker assumption that $X\to Y$ is faithfully flat of finite presentation, with integral fibres. The general case follows from this via Noetherian approximation. Precisely, since $Y$ is quasi-compact and quasi-separated, \SP{01ZA} writes $Y=\lim_i Y_i$ for a cofiltered inverse system $\{Y_i\}$ of finite type integral $\mathbf{Z}$-schemes with affine transition morphisms. Replace each $Y_i$ by its normalization to assume that each $Y_i$ is normal (using \SP{0335}). 
The morphism $\overline{\phi}$ is of finite presentation, and the quasi-compact open immersion $X\hookrightarrow\overline X$ is of finite presentation. Hence, by \SPD{01ZM}{0C0C}, after increasing $i_0$, there exist a smooth morphism $\ov{\phi}_{i_0}\colon \overline{X}_{i_0}\to Y_{i_0}$ and an open subscheme $X_{i_0}\subset\overline X_{i_0}$ whose pullbacks to $Y$ identify with $\overline X$ and $X$. 
By \SP{{0B8W}}, there is an invertible $\sO_{X_{i_0}}$-module $\sL_{i_0}$ whose pullback to $X$ is isomorphic to $\sL$. After increasing $i_0$, the limit argument for morphisms of finite presentation lets us also assume that $\overline{\phi}_{i_0}$ has geometrically integral fibres. For any $i\ge i_0$, denote by $\phi_i\colon X_i\ce X_{i_0}\times_{Y_{i_0}}Y_i\to Y_i$ the base change of $\overline{\phi}_{i_0}|_{X_{i_0}}$ to $Y_i$, and denote by $\sL_i$ the pullback of $\sL_{i_0}$ to $X_i$. By \SPD{01ZM}{01ZP}, any projective embedding of $\overline{X}$ over $Y$ descends to a projective embedding of $\overline{X}_{i}$ over $Y_i$ for large $i$; in particular, $\overline{\phi}_i$ is projective for large $i$.
By limit formalism, for large $i$, $\sL_i$ is trivial when restricted to the generic fibre of $\phi_i$. 
Consequently, for large $i$, the morphism $\phi_i\colon X_i\to Y_i$ and the invertible $\sO_{X_i}$-module $\sL_i$ satisfy all the assumptions of the lemma, so $\sL_i\simeq \phi_i^*\sN_i$ for some invertible $\sO_{Y_i}$-module $\sN_i$. Then $\sL\simeq \phi^*\sN$ where $\sN$ is the pullback of $\sN_i$ to $Y$.
\epf

We are prepared to prove the anticipated Grothendieck--Serre theorem for groups of multiplicative type.
\bthm [cf.~\cite{CTS87}*{4.1--4.3}]\label{G-S type results for mult type}
Let $R$ be a semilocal Pr\"ufer domain, let $A$ be a semilocal domain essentially smooth over $R$, set $X\ce \Spec A$, and write $K(X)=\Frac A$.
For a finite type $X$-group $M$ of multiplicative type, the restrictions
\[
   H^1_{\fppf}(X,M) \hookrightarrow H^1_{\fppf}(K(X),M)\q \text{and}\q  H^2_{\fppf}(X,M) \hookrightarrow H^2_{\fppf}(K(X),M)\q \text{are injective.}
\]
Further, if $M$ is a $\mathrm{flasque}$ $X$-torus, then the following restriction map
 \[
 H^1_{\et}(X,M) \isoto H^1_{\et}(K(X),M) \q \text{ is bijective.}
 \]
\ethm
\bcl\label{conditions}
Let $X$ be a nonempty integral normal scheme essentially smooth over the fixed semilocal Pr\"ufer domain $R$, and let $M$ be a finite type $X$-group of multiplicative type. In the proof of \Cref{G-S type results for mult type}, we replace the assumption $X=\Spec A$ by either of the following conditions.
\benumr 
\item \label{apt-i} There are an essentially smooth semilocal $R$-domain $\Lambda$ and a quasi-compact open immersion $X\hra \ov{X}$ for a smooth projective $\Lambda$-scheme $\ov{X}$ with geometrically integral fibres, such that 
\begin{itemize}
\item $\Pic X_L=0$ for any finite separable field extension $L/\Frac \Lambda$, and
\item $M=N_X$ for a $\Lambda$-group $N$ of multiplicative type.
\end{itemize}
\item\label{apt-ii} There is a connected finite \'etale Galois cover $X\pr\ra X$ splitting $M$ such that any intermediate cover $X\pr\ra X\prpr\ra X$ satisfies
\[
   \Pic X\prpr=0.
\]
\eenum
\ecl
For $X=\Spec A$ as in \Cref{G-S type results for mult type}, \ref{apt-i} is automatic. Indeed, take the trivial compactification $X=\ov{X}$. 
We next check that \ref{apt-i} implies \ref{apt-ii}.
Write $K(-)$ for function fields.
Let $\Lambda\to B$ be a connected finite \'etale Galois cover splitting $N$, put $X\pr\ce X\times_\Lambda B$, and set $\Gamma\ce \Aut_\Lambda(B)$. %
 Since $X_{\Frac \Lambda}$ is geometrically integral, $K(X)\otimes_{\Frac \Lambda}\Frac B$ is a field. Thus $X\pr$ is connected, so $X\pr\to X$ is a connected finite \'etale Galois cover with group $\Gamma$, and it splits $M$.
For an intermediate cover $X\pr\to X\prpr\to X$, put $H\ce \Aut(X\pr/X\prpr)\subset \Gamma$ and $C\ce B^H$.
On function fields,
\[
K(X\prpr)=K(X\pr)^H=(K(X)\otimes_{\Frac \Lambda}\Frac B)^H
       =K(X)\otimes_{\Frac \Lambda}\Frac C=K(X\times_\Lambda C).
\]
The two finite \'etale $X$-schemes $X\prpr$ and $X\times_\Lambda C$ are the normalizations of $X$ in this same field, hence $X\prpr\simeq X\times_\Lambda C$.
For this $C$, the assumption $\Pic X_{\Frac C}=0$ and \Cref{non-Noeth-pullback line bundle}, applied to $X\prpr=X\times_\Lambda C \to \Spec C$, show that the restriction $\Pic C\to \Pic X\prpr$ is surjective.
Since $C$ is semilocal, we conclude that $\Pic C=0=\Pic X\prpr$.
Therefore, \ref{apt-i} falls under \ref{apt-ii}.

\bpf[Proof of \Cref{G-S type results for mult type}]
 Suppose first that $M$ is a flasque $X$-torus. By \Cref{purity for gp of mult type}\ref{gl-H1-H2-flasque}  and a limit argument, the map $H^1_{\et}(X,M) \to H^1_{\et}(K(X),M)$ is surjective.
It remains to prove that the maps $H^i_{\fppf}(X,M)\to H^i_{\fppf}(K(X),M)$ are injective for $i=1,2$.
By \Cref{conditions}, it suffices to work under the assumption \ref{apt-ii}.
We assume first that $M=T$ is an $X$-torus. We choose, as in \cite{CTS87}*{1.3.2} and relative to the splitting cover $X\pr\to X$ in \ref{apt-ii}, a flasque resolution
\[
1\to F \to P \to T \to 1.
\]
Here $F$ is a flasque $X$-torus and $P$ is a quasitrivial $X$-torus, namely a finite product of tori $\mathrm{Res}_{X\prpr/X}\bG_{m,X\prpr}$ with $X\pr \ra X\prpr \to X$ finite \'etale. This yields a commutative diagram
\[
\begin{tikzcd}
{H^1_{\et}(X,P)} \arrow[r] & {H^1_{\et}(X,T)} \arrow[r] \arrow[d, "\rho_1"] & {H^2_{\et}(X,F)} \arrow[d, "\rho_2"] \\
                                   & {H^1_{\et}(K(X),T)} \arrow[r]           & {H^2_{\et}(K(X),F)}
\end{tikzcd}
\]
where the upper row is exact.
By \ref{apt-ii}, we have $H^1_{\et}(X,P)=0$, so the injectivity of $\rho_1$ reduces to that of $\rho_2$.
To prove that $\rho_2$ is injective, we pick $a\in H^2_{\et}(X,F)$ for which $a|_{K(X)}=0$. By spreading out, we may assume that $X$ is a localization of an irreducible, smooth, affine $R$-scheme $\wt{X}$, $F=\wt{F}_X$ for a flasque $\wt{X}$-torus $\wt{F}$, and $a=\wt{a}|_X$ for some class $\wt{a} \in H^2_{\et}(\wt{X},\wt{F})$. 
Since $\wt{a}|_{K(X)}=0$, there exists a proper closed subset $Y\subset \wt{X}$ such that
$$
\wt{a}|_{\wt{X}\backslash Y}=0 \in H^2_{\et}(\wt{X}\backslash Y,\wt{F}).
$$
By \Cref{purity for gp of mult type}\ref{gl-H1-H2-flasque}, $\widetilde{a}=0$, so $a=\widetilde{a}|_X=0$. This proves the injectivity of $\rho_2$ and hence also of $\rho_1$. 

Now let $M$ be an arbitrary $X$-group of multiplicative type. There is an $X$-subtorus $T\subset M$ such that $\mu\ce M/T$ is $X$-finite. Consequently, for any generically trivial $M$-torsor $\cP$, the $\mu$-torsor $\cP/T$ is finite over $X$. Since $X$ is normal, one has $(\cP/T)(X)=(\cP/T)(K(X))$. Therefore, $\cP/T\to X$ has a section that lifts to a generic section of $\cP\to X$, that is, $\cP$ reduces to a generically trivial $T$-torsor $\cP_T$. By the injectivity of $\rho_1$,  $\cP_T$ is trivial, and hence so is $\cP$. This proves the injectivity of $H^1_{\fppf}(X,M) \hra H^1_{\fppf}(K(X),M)$.

Finally, we prove that $H^2_{\fppf}(X,M)\ra H^2_{\fppf}(K(X),M)$ is injective. 
By \cite{CTS87}*{1.3.2}, there is a short exact sequence
\[
1\to M\to F \to P \to 1
\]
of $X$-groups of multiplicative type, where $F$ is flasque and $P$ is quasitrivial, both split after base change by $X\pr\to X$. This yields the following commutative diagram whose upper row is exact
\[
\begin{tikzcd}
{H^1_{\fppf}(X,P)} \arrow[r] & {H^2_{\fppf}(X,M)} \arrow[r] \arrow[d,"\rho_3"] & {H^2_{\fppf}(X,F)} \arrow[d,"\rho_2"] \\
                                   & {H^2_{\fppf}(K(X),M)} \arrow[r]           & {H^2_{\fppf}(K(X),F).}
\end{tikzcd}
\]
 Since we have already shown that $H^1_{\fppf}(X,P)=0$ and $\rho_2$ is injective by the torus case treated above, which used \Cref{purity for gp of mult type}\ref{gl-H1-H2-flasque}, the injectivity of $\rho_3$ follows.
\epf

\section{Geometric lemmata}
\label{sect-geom lem}
\bpp[Geometric presentation lemma in the style of Gabber--Quillen]
We begin with the following Pr\"uferian analogue of Gabber--Quillen type results, which play a prominent role in both Fedorov's and $\check{\mathrm{C}}$esnavi$\check{\mathrm{c}}$ius' work on mixed characteristic Grothendieck--Serre, see \cite{Fed22b}*{Proposition~3.18} and \cite{Ces22a}*{Variant~3.7}, respectively.
\epp

\blemt \label{Ces's Variant 3.7}
Let $R$ be a semilocal Pr\"ufer domain. Fix the following data:
\begin{itemize}
\item [--] a projective, flat $R$-scheme $\overline{X}$ with fibres of pure dimension $d>0$,
\item [--] an open subscheme $X\subset \overline{X}$ smooth over $R$, a finite subset $\mathbf{x}\subset X$, and
\item [--] a closed subscheme $Y\subset X$ and its closure $\ov{Y}$ in $\ov{X}$.
\end{itemize}
If $Y$ is $R$-fibrewise of codimension $\ge 1$ in $X$ and $\overline{Y}\backslash Y$ is $R$-fibrewise of codimension $\ge 2$ in $\overline{X}$, there are
\begin{itemize}
\item an affine open $S\subset \mathbb{A}_R^{d-1}$ and an affine open $U\subset X$ containing $\mathbf{x}$, and
\item a smooth morphism $\pi\colon U\to S$ of relative dimension 1 such that $Y \cap U$ is $S$-finite.
\end{itemize}
\elemt
\bpf
We first reduce to a finite dimensional base.  By \Cref{geom}\ref{limitarg} and the limit arguments \cite{EGAIV3}*{th\'eor\`eme~8.8.2, proposition~11.2.6, and proposition~9.2.6.1}, the data spread out to a semilocal Pr\"ufer subdomain of finite Krull dimension.  Enlarging this subdomain if needed, we may arrange that the descended $\overline X$ is projective and flat with fibres of pure dimension $d$, that the descended $X$ is smooth, and that the two fibrewise codimension conditions hold.  The finite set $\mathbf{x}$ also descends.  Thus we may assume that $R$ has finite Krull dimension.

Following \cite{Ces22a}*{proof of Variant~3.7}, we reduce to the case where $\mathbf{x}$ lies over closed points of $\Spec R$.  
First replace each $x\in\mathbf{x}$ by a specialization $x'$ closed in its $R$-fibre.  Every open neighbourhood of $x'$ also contains $x$.  Thus it remains to treat one $x$ over a nonmaximal $\fp\subset R$.

Write $R_{\fp}=\bigcup_A A$ with $A\subset R_{\fp}$ finite type over $\mathbb{Z}$, and put $\fp_A=A\cap\fp$.  By \SPD{0EY1}{0C0C} and \cite{EGAIV3}*{proposition~9.2.6.1}, after enlarging $A$, the $R_{\fp}$-data, together with an ample line bundle on $\overline X_{R_{\fp}}$, descend over $A$, preserving flatness, smoothness, pure dimensionality, and the codimension bounds. 
Also, $x$ descends to an $A/\fp_A$-finite closed subscheme $\widetilde{x}$ of the descended $X$.

Enlarge $A$ so that $(\fq\backslash\fp)\cap A^\times\neq\emptyset$ for every $\fq\supset\fp$ with $\operatorname{ht}(\fq)=\operatorname{ht}(\fp)+1$.  Then $A\cdot R_{\fm}=R_{\fp}$ for every maximal $\fm\supset\fp$.  Since $R_{\fm}/\fp R_{\fm}$ is a nontrivial valuation ring of $k_{\fp}$, the field $k_{\fp}$ is not algebraic over a finite field.  From $k_{\fp}=\bigcup_A A/\fp_A$, after enlarging $A$ again, $A/\fp_A$ is not algebraic over a finite field.  Hence $A/\fp_A$ has a nonzero prime $\fp\pr$\footnote{We use the following fact: a prime ideal of a finite type $\mathbb{Z}$-algebra is maximal if and only if its residue field is finite.}.

Using \cite{Bou98}*{chapitre~VI, \S3, \textnumero~1, th\'eor\`eme~1}, choose a valuation ring $V_{\fp\pr}\subset k_{\fp}$ whose restriction to $k_{\fp_A}\ce\Frac(A/\fp_A)$ has centre $\fp\pr$.  Let $V\subset\Frac R$ be its preimage under $R_{\fp}\twoheadrightarrow k_{\fp}$.  Then $V\cdot R_{\fq}=R_{\fp}$ for every $\fq\supset\fp$.  By \cite{Bou98}*{chapitre~VI, \S7, \textnumero~1, propositions~1~et~2},
\[
        R\pr\ce R\cap V
\]
is a semilocal Pr\"ufer domain and $\Spec R\pr$ is obtained by gluing $\Spec R$ and $\Spec V$ along $\Spec R_{\fp}$.

We base change the descended $(\overline X,X,Y,\overline Y)$-data from $A$ to $V$, and glue them with the original $(\overline X,X,Y,\overline Y)$-data over $R$ along their common $R_{\fp}$-base change.  The resulting $R\pr$-data still have $\overline X$ projective and flat with fibres of pure dimension $d$, with $X$ smooth and the same fibrewise codimension bounds.  Since $\widetilde{x}_V$ is $V$-finite, $x$ specializes to a point over the closed point of $\Spec V$.  Proving the lemma after this replacement suffices, because restriction to $\Spec R\subset\Spec R\pr$ recovers an open containing the original $x$.  Iterating over $\mathbf{x}$ proves the reduction.  The same argument with affine smooth $R$-schemes in place of projective compactifications gives the affine version used in \Cref{reduction to closed points}.  The old closed fibres are unchanged and each new closed fibre contains only the specialization chosen at that step, so the cardinality assumption in \Cref{variant of Lindel's lem} is preserved.

Choose a projective embedding $\overline X\hookrightarrow \mathbb P_R^N$.  
After replacing $\sO_{\overline X}(1)$ by a large power, fix positive weights $w_1,\ldots,w_{d-1}$ so that, on every closed fibre, the weighted tuples satisfying \cite{Ces22a}*{Proposition~3.6(i)--(vi)} form a nonempty open.  We also arrange that ambient forms of degrees $1,w_1,\ldots,w_{d-1}$ restrict surjectively to the corresponding spaces of sections on every closed fibre.  Let $\mathcal P$ be the affine $R$-space of tuples of ambient homogeneous forms $h_0,\ldots,h_{d-1}$ of degrees $1,w_1,\ldots,w_{d-1}$.  Such a tuple defines the weighted rational map
\[
        \bar p\colon \overline X\dashrightarrow \mathbb P_R(1,w_1,\ldots,w_{d-1}).
\]
On the chart $h_0\ne 0$ it is a map to $\mathbb A_R^{d-1}$.  For a maximal ideal $\fm\subset R$, let $\Omega_{\fm}\subset \mathcal P_{k_{\fm}}$ be the open locus where \cite{Ces22a}*{Proposition~3.6(i)--(vi)} holds for the triple $(\overline X_{k_{\fm}},X_{k_{\fm}},\overline Y_{k_{\fm}})$ and for the finite set $\mathbf{x}_{k_{\fm}}$.  The hypotheses of that proposition are satisfied because $Y$ is fibrewise of codimension $\ge 1$ and $\overline Y\backslash Y$ is fibrewise of codimension $\ge 2$.  The weights have been chosen, as in the proof of the cited proposition and uniformly over the finitely many closed fibres, so that $\Omega_{\fm}(k_{\fm})\ne\emptyset$ for all $\fm$.  When $k_{\fm}$ is infinite, a nonempty open of an affine space has a $k_{\fm}$-point.  When $k_{\fm}$ is finite, the finite-field Bertini input in the cited proof supplies such a point after increasing the weights.
Choose $\omega_{\fm}\in\Omega_{\fm}(k_{\fm})$ for every $\fm$.  Since $R$ is semilocal and $\mathcal P$ is an affine space, the Chinese remainder theorem lifts $(\omega_{\fm})_{\fm}$ to an $R$-point of $\mathcal P$.  We keep the corresponding tuple $h_0,\ldots,h_{d-1}$ and the associated map $\bar p$.

The conditions checked on closed fibres hold over $R$, because any nonempty closed subset of the proper $R$-scheme $\overline X$ meets a closed fibre.  Thus the base locus of $\bar p$ is disjoint from $\mathbf{x}\cup\overline Y$, and
\[
(\overline Y\backslash Y)\cap \bar p^{-1}\bigl(\bar p(\mathbf{x})\bigr)=\emptyset
\quad\text{and}\quad
\overline Y\cap V(h_0)\cap \bar p^{-1}\bigl(\bar p(\mathbf{x})\bigr)=\emptyset.
\]
By \cite{Ces22a}*{Proposition~3.6(iii), (v)} and openness of the smooth locus, $\bar p$ is smooth of relative dimension $1$ at the points of $\mathbf{x}\cup\bigl(\overline Y\cap \bar p^{-1}(\bar p(\mathbf{x}))\bigr)$.
Let $T\subset \mathbb A_R^{d-1}$ be the finite set $\bar p(\mathbf{x})$.  The two displayed equalities mean that $T$ is disjoint from the closed images of $\overline Y\backslash Y$ and $\overline Y\cap V(h_0)$ under $\bar p$. Thus, after replacing $\mathbb A_R^{d-1}$ by a principal open neighbourhood $S$ of $T$, we have
\[
(\overline Y\backslash Y)\cap \bar p^{-1}(S)=\emptyset
\quad\text{and}\quad
\overline Y\cap V(h_0)\cap \bar p^{-1}(S)=\emptyset.
\]
By the argument proving \cite{Ces22a}*{Proposition~3.6(vii)}, after shrinking $S$ around $T$, the morphism $\overline Y\cap \bar p^{-1}(S)\to S$ is quasi-finite.  It is proper, hence finite by \SP{02OG}.  We shrink $S$ once more so that $\bar p$ is smooth of relative dimension $1$ along this finite $S$-scheme and along $\mathbf{x}$.

Choose an affine open neighbourhood $U\subset X\cap \bar p^{-1}(S)$ of $\mathbf{x}\cup(\overline Y\cap \bar p^{-1}(S))$ inside the smooth locus of $\bar p$.  This is possible inside the affine scheme $(\overline X\backslash V(h_0))\cap \bar p^{-1}(S)$.  Then $\pi\ce\bar p|_U\colon U\to S$ is smooth of relative dimension $1$.  Since $(\overline Y\backslash Y)\cap \bar p^{-1}(S)=\emptyset$, we have $Y\cap U=\overline Y\cap \bar p^{-1}(S)$, so $Y\cap U$ is $S$-finite.
\epf

\bpp[A variant of Lindel's lemma]
   A useful lemma due to Lindel \cite{Lin81}*{Proposition~1 and Lemma} states that an \'etale extension of local rings $A\to B$ with trivial extension of residue fields induces isomorphisms along a well-chosen non-unit $r\in A$:
   $$
   A/r^nA\isoto B/r^nB, \q \text{ where } \q n\ge1.
   $$
	   In our context, where $B$ is essentially smooth over a valuation ring, we will prove a variant of \emph{loc.~cit.} that allows us to fix the $r\in B$ in advance, at the cost of carefully choosing $A$ to be a local ring of an affine space over that valuation ring. This result will be a crucial geometric tool for addressing torsors under a `constant' reductive group scheme. Like Lindel's work on the Bass--Quillen conjecture for vector bundles and more generally our work \cite{GL25} on the Bass--Quillen conjecture for torsors, it also reduces our problem to studying torsors on open subsets of relative affine spaces.
\epp

\bpropt \label{variant of Lindel's lem}
Let $\Lambda$ be a semilocal Pr\"ufer domain. We fix
\begin{itemize}
\item [--]an irreducible, affine, smooth $\Lambda$-scheme $X$ of pure relative dimension $d>0$, 
\item [--]a finitely presented closed subscheme $Y\subset X$ that is $\Lambda$-fibrewise of codimension $>0$, and 
\item [--]a finite subset $\mathbf{x} \subset X$.
\end{itemize}
If each fibre of $\mathbf{x}$ over a maximal ideal $\fm\subset \Lambda$ has fewer than $\max(\#\,k_{\mathfrak{m}},d)$ points, then there are 
\begin{itemize}
\item an affine open $W\subset X$ containing $\mathbf{x}$,
\item an affine open $U\subset \mathbb{A}_{\Lambda}^d$, and
\item an \'etale $\Lambda$-morphism $f\colon W\to U$ fitting into the following Cartesian square 
\begin{equation*}
\qqq\begin{tikzcd}
W\cap Y \arrow[r, hook]  \arrow[d, equal]
 & W \arrow[d, "f"] \\
W\cap Y \arrow[r, hook]
& U, 
\end{tikzcd}
\x{\qqq where $f|_{W\cap Y}\colon W\cap Y \to U$ is a closed immersion.}
\end{equation*}
\end{itemize}
\epropt

\bremt\label{rem on Lindel's lem}
The assumption on $\# \mathbf{x}$ holds, for instance, if either $\mathbf{x}$ is a singleton or $d>\# \mathbf{x}$. 
The latter case $d>\# \mathbf{x}$ is crucial for the general semilocal setup of \Cref{G-S for constant reductive gps}.
Note that a certain assumption on $\# \mathbf{x}$ is necessary: when $X$ is a smooth affine curve over a finite field $\GGL=\bF_q$ and $\mathbf{x}\subset X(\bF_q)$, the resulting map $f$ from \Cref{variant of Lindel's lem} would force an injection $\mathbf{x} \hookrightarrow \mathbb{A}_{\bF_q}^1$, which is impossible as soon as $\# \mathbf{x} > q$.
\eremt
To prove \Cref{variant of Lindel's lem}, we begin with the following reduction to considering closed points.
\blemt \label{reduction to closed points}
It suffices to prove \Cref{variant of Lindel's lem} under the following supplementary assumptions:
\benumr
\item\label{red-ass-i} $Y=\{a=0\}$ for a function $a\in \GG(X,\sO_X)$ such that $\mathbf{x}\subset Y$, and
\item\label{red-ass-ii} every $x\in \mathbf{x}$ is a closed point in some closed $\Lambda$-fibre of $X$.
\eenum
\elemt
\bpf
Assume that \ref{red-ass-ii} holds. 
As the relative dimension of $X$ is nonzero, the subset $\mathbf{x}\cup Y$ avoids all the maximal points of the $\Lambda$-fibres of $X$.
     By prime avoidance, there is an $a\in \Gamma(X,\sO_X)$ whose vanishing locus $V(a)$ contains $\mathbf{x}\cup Y$ and has nowhere dense $\GGL$-fibres in $X$.
    Replacing $Y$ by $V(a)$ achieves the reduction to \ref{red-ass-i}.

It remains to reduce to the assumption \ref{red-ass-ii}.
First, using a limit argument (\Cref{geom}\ref{limitarg}), we reduce to the case where $\GGL$ has finite Krull dimension.  Suppose that, for each $x\in \mathbf{x}$, the closure $\ov{\{x\}}$ contains a closed point $z_x$ in a closed $\GGL$-fibre of $X$.  Put $\mathbf{z}\ce\{z_x\}_{x\in\mathbf{x}}$, and assume that $\mathbf{z}$ satisfies the same cardinality assumption.  Then it is enough to prove the proposition for $\mathbf{z}$, because every open neighbourhood of $z_x$ contains the generalisation $x$.

When such specializations are not available, the affine form of the auxiliary-prime step in the proof of \Cref{Ces's Variant 3.7} supplies them after replacing the base by a semilocal Pr\"ufer domain whose spectrum contains the old one as an open, and it preserves the cardinality assumption.  Therefore the preceding paragraph reduces us to \ref{red-ass-ii}.
\epf

\blemt \label{lem on g}
In the situation of \Cref{reduction to closed points}, there exists a $\GGL$-morphism $g\colon X\to \mathbb{A}_{\GGL}^{d-1}$ such that
\benumr
\item\label{redc-i} it is \emph{smooth of relative dimension $1$} at $\mathbf{x}$.
\item\label{redc-ii} the restriction $g|_Y$ is \emph{quasi-finite} at $\mathbf{x}$.
\eenum
In addition, if $d>\#(\mathbf{x} \cap X_{k_{\mathfrak{m}}})$ for every maximal ideal $\fm\subset \Lambda$ with finite residue field, then we may choose $g$ so that, for every such maximal ideal $\fm$, the points of $\mathbf{x}\cap X_{k_{\mathfrak{m}}}$ have pairwise distinct $k_{\mathfrak{m}}$-rational images.
\elemt

\bpf
We first reduce the lemma to the case when $\GGL=k$ is a field. Suppose that for every maximal ideal $\mathfrak{m}\subset \GGL$, there exists a $k_{\mathfrak{m}}$-morphism $g_{\mathfrak{m}}\colon X_{k_{\mathfrak{m}}}\to \mathbb{A}_{k_{\mathfrak{m}}}^{d-1}$ that is smooth at $\mathbf{x}\cap X_{k_{\mathfrak{m}}}$ such that the restriction $g_{\mathfrak{m}}|_{Y_{k_{\mathfrak{m}}}}$ is quasi-finite at $\mathbf{x}\cap X_{k_{\mathfrak{m}}}$, and, when $k_{\mathfrak m}$ is finite and $d>\#(\mathbf{x}\cap X_{k_{\mathfrak m}})$, sends the points of $\mathbf{x}\cap X_{k_{\mathfrak m}}$ to pairwise distinct $k_{\mathfrak m}$-rational points.
We then use the Chinese remainder theorem to lift the maps $\{g_{\mathfrak{m}}\}_{\mathfrak{m}}$ simultaneously to a $\GGL$-morphism $g\colon X\to \mathbb{A}_{\GGL}^{d-1}$ satisfying \ref{redc-i}--\ref{redc-ii} (note that the smoothness of $g$ at $\mathbf{x}$ follows from the same property for its closed fibres and the fibral criterion for flatness \cite{EGAIV3}*{théorème~11.3.10}).  For such $\mathfrak m$, this lift restricts to the chosen $g_{\mathfrak m}$, so the images remain pairwise distinct and $k_{\mathfrak m}$-rational.

For the rest, we begin with the case when the field $\Lambda=k$ is perfect.
By \Cref{reduction to closed points}~\ref{red-ass-ii}, every element of the finite set $\mathbf{x}$ is a closed point of $X$.  We replace $\mathbf{x}$ by the reduced closed subscheme $\coprod_{x\in\mathbf{x}}\Spec k_x\subset X$, where $k_x=\kappa(x)$, and keep the same notation. Choose a map of sets
\[
  \tau=(\tau_1,\ldots,\tau_{d-1})\colon |\mathbf{x}|\to k^{d-1}.
\]
For each $j$, the tuple $(\tau_j(x))_{x\in|\mathbf{x}|}\in\prod_{x\in|\mathbf{x}|}k_x=\Gamma(\mathbf{x},\sO_{\mathbf{x}})$ defines a coordinate function $t_j$ and hence a $k$-morphism $t=(t_1,\ldots,t_{d-1})\colon \mathbf{x}\to \bA^{d-1}_k$.
For $1\le i\le d-1$, let $T_i\subset \bA_k^i$ be the finite reduced image of $(t_1,\ldots,t_i)\colon \mathbf{x}\to \bA_k^i$.  Put $T_0=\Spec k$.
We will construct $g$ inductively with $g|_{\mathbf{x}}=t$.
If $k$ is finite, choose $\tau$ injective, as allowed by our cardinality assumption.
\bclt \label{inductively construct g_i}
	 For every $1\le i\le d-1$, there exists a $k$-morphism $g_i\colon X\to \bA_k^i$ such that
\begin{itemize}
  \item $g_i$ is \emph{smooth} at $\mathbf{x}$ and $g_i|_{\mathbf{x}}=(t_1,\ldots,t_i)$.
  \item  $(g_i|_Y)^{-1}(T_i)$ has \emph{pure dimension $d-1-i$}.
\end{itemize}
\eclt
\bpf [Proof of the claim]
Let $g_0\colon X\to \Spec k$ be the structural morphism. 
Assume by induction that $g_{i-1}\colon X\ra \bA^{i-1}_k$ has been constructed.
Put $Y_{i-1}\ce (g_{i-1}|_Y)^{-1}(T_{i-1})$.
For $x\in\mathbf{x}$, the point $g_{i-1}(x)$ is $k$-rational by induction, and the fibre $g_{i-1}^{-1}(g_{i-1}(x))$ is smooth over $k$ at $x$.  As $k$ is perfect, $k_x/k$ is separable, and the completed local ring of this fibre at $x$ is a power series ring $k_x\fps{u_1,\ldots,u_{d-i+1}}$ by \SP{07EJ}.
Let $x^{(1)}$ be the closed subscheme defined by $\fm_x^2$ in this fibre, where $\fm_x$ is the maximal ideal of its local ring at $x$.  We may choose a section $a_x\in \Gamma(x^{(1)},\sO_{x^{(1)}})$ whose residue is $t_i(x)\in k_x$ and such that the class of $a_x-t_i(x)$ in $\fm_x/\fm_x^2$ is nonzero.
Because the $x^{(1)}$ are disjoint finite closed subschemes of the affine scheme $X$, the Chinese remainder theorem produces an $h_0\in \Gamma(X,\sO_X)$ whose restriction to each $x^{(1)}$ is $a_x$.

It remains to impose the dimension condition on $Y_{i-1}$.
Let $\eta_1,\ldots,\eta_N$ be the generic points of the irreducible components of $Y_{i-1}$.
Let $I\subset\Gamma(X,\sO_X)$ be the ideal defining $\bigsqcup_{x\in\mathbf{x}}x^{(1)}$.  By prime avoidance, choose $h\in h_0+I$ such that, for every $\eta_j$ and every value $\lambda\in t_i(|\mathbf{x}|)\subset k$, the function $h-\lambda$ is nonzero at $\eta_j$.
Since $h|_{x^{(1)}}=a_x$ for every $x\in\mathbf{x}$, the restriction of $h$ to the fibre $g_{i-1}^{-1}(g_{i-1}(x))$ is smooth at $x$, and hence $(g_{i-1},h)$ is smooth at every $x\in\mathbf{x}$.
Then no irreducible component of $Y_{i-1}$ is contained in $h^{-1}(t_i(\mathbf{x}))$, so Krull's principal ideal theorem shows that $Y_{i-1}\cap h^{-1}(t_i(\mathbf{x}))$ has pure dimension $d-1-i$.
Taking $g_i\ce (g_{i-1},h)\colon X\to \bA_k^i=\bA_k^{i-1}\times_k\bA_k^1$ completes the induction.
\epf

Let $g\ce g_{d-1}$ be the morphism of \Cref{inductively construct g_i}. The smoothness assertion in the claim proves \ref{redc-i}. Since $g|_{\mathbf{x}}=t$, the set $T_{d-1}$ is $g(\mathbf{x})$. The claim also says that $(g|_Y)^{-1}(T_{d-1})$ has pure dimension $0$. As $g|_Y$ is of finite type, it is quasi-finite at the points of $\mathbf{x}$, proving \ref{redc-ii}. If $k$ is finite and $\tau$ is injective, then for every $x\in|\mathbf{x}|$ the point $g(x)=t(x)$ is the $k$-point $\tau(x)$ of $\mathbb A_k^{d-1}$. Hence the images of the points of $\mathbf{x}$ are pairwise distinct $k$-points.

It remains to consider an imperfect field $k$. Then $k$ is infinite, so the pairwise distinctness assertion, which is only imposed over finite residue fields, does not apply. The required morphism follows from the geometric presentation of Gabber--Gros--Suwa \cite{CTHK97}*{Theorem 3.1.1}, which even ensures that $g|_Y$ is finite.
\epf

We can now finish the proof of \Cref{variant of Lindel's lem}.
\bpf
[Proof of \Cref{variant of Lindel's lem}]
 Let $g\colon X\ra \bA_{\Lambda}^{d-1}$ be as in \Cref{lem on g}.
By the openness of the quasi-finite locus (resp., smooth locus) of a finitely presented morphism, we can shrink $X$ around $\mathbf{x}$ to assume that $g$ is smooth of relative dimension $1$ and that $g|_Y$ is quasi-finite. 
If $\fm \subset \GGL$ is  a maximal ideal with finite residue field, then by the assumption of \Cref{variant of Lindel's lem}, either $\#(\mathbf{x}\cap X_{k_{\fm}})\le \#\,k_{\fm}-1$, or every fibre of $g_{k_{\fm}}$ contains at most 1 point of $\mathbf{x}$ by the additional assertion in \Cref{lem on g}. 
Consequently, for every closed point $s$ of $\Spec\sO_{\bA_{\Lambda}^{d-1},g(\mathbf{x})}$, the finite reduced fibre $\mathbf{x}\times_{\bA_{\Lambda}^{d-1}}s$ admits a closed immersion into $\bA^1_s$ whose complement has a $k(s)$-point.  For infinite $k(s)$, choose a generator of the finite reduced $k(s)$-algebra while avoiding one further $k(s)$-value.  For finite $k(s)$, the previous sentence leaves such a value unused. Using a limit argument, \Cref{variant of Lindel's lem} then follows by applying the following \Cref{quasifinite embedding} to the base change of $g$ to the semilocal ring $R\ce \sO_{\bA_{\Lambda}^{d-1},g(\mathbf{x})}$, with $Z\ce \mathbf{x}$.
  \epf
  
  As pointed out by Fedorov in \cite{Fed23}*{Proof of Proposition 4.5}, our original argument for \Cref{variant of Lindel's lem} actually leads to the following result (cf.~\cite{CF23}*{Lemma 2.5}).
  
\begin{variant}\label{quasifinite embedding}
Let $R$ be a semilocal ring with spectrum $S$, and suppose that we are given
\begin{itemize}
    \item[--] a smooth affine morphism $C \to S$ of pure relative dimension 1, and
    \item[--] closed subschemes $Z\subset Y \subset C$ such that $Z$ is $S$-finite and $Y$ is $S$-quasi-finite.
\end{itemize}
If, for each closed point $s\in S$, there is an embedding $Z_s \hra \mathbb{A}^1_{s}$ such that $(\mathbb{A}^1_{s}\backslash Z_{s})(k_s)\neq \emptyset$, then there are 
\begin{itemize}
    \item affine opens $U\subset \bA^1_S$ and $W \subset C$, with $W$ containing $Z$, and
    \item an \'etale $S$-morphism $f\colon W\to U$ fitting into the following Cartesian square 
\begin{equation*}
\qqq\begin{tikzcd}
W\cap Y \arrow[r, hook]  \arrow[d, equal]
 & W \arrow[d, "f"] \\
W\cap Y \arrow[r, hook]
& U, 
\end{tikzcd}
\x{\qqq where $f|_{W\cap Y}\colon W\cap Y \to U$ is a closed immersion.}
\end{equation*}
\end{itemize}
\end{variant}

\bpf
By prime avoidance for the ideal of $Y$ and the maximal points of the closed fibres of $C/S$, choose an $a\in \Gamma(C,\sO_C)$ vanishing on $Y$ such that $\{a=0\}$ has finite closed $S$-fibres. View $a$ as an $S$-morphism $C\to \bA_S^1$. Its quasi-finite locus is open \cite{EGAIV4}*{théorème~13.1.3}. For a point of $Z$ over a closed $s\in S$, the fibre of $a$ over $(s,0)$ is $\{a=0\}_s$, hence finite. Since $Z$ is $S$-finite and $S$ is semilocal, the remaining points of $Z$ are generizations of such points. Thus the quasi-finite locus contains $Z$. Shrinking $C$ to an open neighbourhood of $Z$ inside this locus, $a$ becomes quasi-finite, and so $\{a=0\}$ is $S$-quasi-finite. Replacing $Y$ by $\{a=0\}$, we may assume that $Y=\{a=0\}$.

By Zariski's main theorem \cite{EGAIV4}*{corollaire~18.12.13}, $a$ factorises as an open immersion $j$ followed by a finite morphism $h_1$:
\[
  a=h_1\circ j\colon C \overset{j}{\hra} \overline{C} \xrightarrow {h_1} \bA_S^1. 
\]
We regard $h_1$ as an element of $\Gamma(\ov{C},\sO_{\ov{C}})$.
Write $S_0\subset S$ for the union of the closed points of $S$ (with the reduced structure). 
Let $(-)_0$ denote the base change functor from $S$ to $S_0$.

\bclt \label{lem on b}
   There exists an $h_2\in \Gamma(\ov{C}, \sO_{\ov{C}})$, regarded as a morphism $ \overline{C} \to \bA_S^1$, satisfying the following:
\benumr
\item \label{lem on b 1}  we have $h_1^{-1}(h_1(Z_0)) \cap h_2^{-1}(h_2(Z_0))= Z_0$ set-theoretically.
\item \label{lem on b 2} $h_2$ is \'etale around $Z$ and induces a bijection $Z_0 \xrightarrow{\sim} h_2(Z_0)$ of sets.
\item \label{lem on b 3} $h_2$ induces an isomorphism of residue fields $k_{h_2(x)} \xrightarrow{\sim} k_x$ for every $x \in Z_0$.
\eenum
\eclt
\bpf[Proof of the claim] 
  For $x\in Z_0$, let $x^{(1)}$ be the first infinitesimal neighbourhood of $x$ in the smooth fibre $C_s$, where $s\in S_0$ is the image of $x$, and set $Z_0^{(1)}\ce \bigsqcup_{x\in Z_0}x^{(1)}$.
  Applying \cite{CF23}*{Lemma~2.6} on each fibre over $S_0$ to the given closed immersions $Z_s\hra\bA^1_s$, we obtain a closed immersion of $S_0$-schemes
\begin{equation} \label{embeds first neighb of x}
 \tst \iota_1\colon Z_0^{(1)} \hra \bA_{S_0}^1 \subset \bA_S^1
\end{equation}
whose restriction to $Z_0$ is the given closed immersion. In particular, it induces a bijection onto its image and the required isomorphisms of residue fields.
Moreover, by the assumption of \Cref{quasifinite embedding}, the complement of $\iota_1(Z_s)$ has at least one $k_s$-rational point for each $s\in S_0$.
Set
\[
  E\ce \left(h_1^{-1}(h_1(Z_0)) \backslash Z_0 \right)_{\mathrm{red}}.
\]
For each $s\in S_0$, choose a $k_s$-rational point of $\bA^1_s\backslash \iota_1(Z_s)$, and map $E_s$ constantly to this point. This gives a section on the finite reduced scheme $E$ over $S_0$ whose values avoid $\iota_1(Z_0)$.
Together with \eqref{embeds first neighb of x}, this determines a section $\iota\in\Gamma(T,\sO_T)$ on the finite closed subscheme
\[
  \tst T\ce Z_0^{(1)}\sqcup E\subset \overline C
\]
such that $h_1^{-1}(h_1(Z_0))\cap \iota^{-1}(\iota(Z_0))=Z_0$ set-theoretically.

Let $h_2\in \Gamma(\ov{C},\sO_{\ov{C}})$ be a lifting of $\iota$.  Such a lifting exists because $\overline C$ is affine over $S$. 
The construction proves \ref{lem on b 1}, \ref{lem on b 3}, and the bijectivity assertion in \ref{lem on b 2}. 
The restriction of $h_2$ to the first infinitesimal neighbourhoods agrees with $\iota_1$, so the maps on the fibres over closed points of $S$ are \'etale at the points of $Z_0$.  Thus $h_2$ is \'etale at these points by the fibre criterion for \'etaleness. 
Since the \'etale locus is open and every point of the $S$-finite scheme $Z$ is a generization of a point of $Z_0$, this open locus contains $Z$.
 \epf
 
 Let $B\ce R[h_1,h_2]\subset \Gamma(\ov{C},\sO_{\ov{C}})$.
 Set $V\ce \Spec B$. 
For $i=1,2$, let $v_i\colon R[t]\ra B$ be the map $t\mapsto h_i$.
It induces the factorization $h_i=\Spec(v_i)\circ h_3$ for $h_3\colon \ov{C}\to V$.
Denote $Y\pr \ce {\Spec(B/(h_1))} \subset V$. 
As $h_1|_C=a$, we have the following equation
\begin{equation}\label{eq-Y}
   ( h_3|_C)^{-1}(Y\pr)=\{a=0\}=Y.
\end{equation}
The morphism $h_3$ is finite. Indeed, $\Gamma(\overline C,\sO_{\overline C})$ is finite over $R[h_1]$, so $B$ is finite over $R[h_1]$ and $\Gamma(\overline C,\sO_{\overline C})$ is finite over $B$.
Note that $v_2$ restricts to a closed immersion $\ov{v}_2\colon Y\pr \hra \bA_R^1 $ since we have $R[h_2]\twoheadrightarrow B/(h_1)$.

 \bclt \label{isom of loc rings}
 The base change of $h_3$ to $\sO_{V,h_3(Z_0)}$ is an isomorphism of semilocal rings \[\sO_{V,h_3(Z_0)} \xrightarrow{\sim} \sO_{\ov{C},Z_0}\,=\sO_{C,Z_0}.\]
 \eclt
 
\bpf[Proof of the claim]
The map is injective by the definition of $B\subset \Gamma(\overline C,\sO_{\overline C})$, so it remains to prove surjectivity. 
The \Cref{lem on b}~\ref{lem on b 2}--\ref{lem on b 3} imply that $h_3$ induces a bijection $Z_0 \isoto h_3(Z_0)$ of sets and an isomorphism of residue fields $k_{h_3(x)} \isoto k_x$ for every $x\in Z_0$. Moreover, \Cref{lem on b}~\ref{lem on b 1}--\ref{lem on b 2} imply that
$Z_0 = h_3^{-1}(h_3(Z_0)) $ as sets. The morphism $h_3$ is finite because $h_1$ is finite, and it is unramified along $Z_0$ because $h_2$ is \'etale there. Thus Nakayama's lemma applies. Indeed, if $J$ is the Jacobson radical of $\sO_{V,h_3(Z_0)}$, then the following composite
\[
\text{$\tst \prod_{x\in Z_0}k_{h_3(x)} \simeq \sO_{V,h_3(Z_0)}/J \xrightarrow{h_3^*} \sO_{\ov{C},Z_0}/J\sO_{\ov{C},Z_0} \simeq \prod_{x\in Z_0} k_x$ \qq   is bijective}. \qedhere
\]
\epf

Shrinking $C$ to the \'etale locus of $h_2$ around $Z$, and replacing $Y$ by its intersection with this open, we may assume that
\[
\phi\ce h_2|_C\colon C\to \mathbb{A}_{R}^1
\]
is \'etale.
By \Cref{isom of loc rings}, there exists an affine open neighbourhood $W\pr \subset V$ of $h_3(Z_0)$ such that a restriction of $h_3$ induces $h_3^{-1}(W\pr ) \isoto  W\pr$. 
As noted above, $v_2\colon V\ra \bA^1_R$ restricts to a closed immersion $\ov{v}_2 \colon Y\pr \hra  \mathbb{A}_{R}^1$.
 In particular, the topology of $Y\pr $ is induced from that of $\mathbb{A}_{R}^1$ via $\ov{v}_2$. 
 Consequently, there exists an affine open neighbourhood $U\subset \bA_R^1$ of $\phi(Z_0)$ such that $\ov{v}_2^{-1}(U)\subset W\pr $. 
 Then $\phi$ induces a closed immersion of affine schemes
 \[
 Y_U\ce \phi^{-1}(U)\cap Y=(v_2\circ h_3|_C)^{-1}(U)\cap Y \overset{(\ref{eq-Y})}{\simeq}  v_2^{-1}(U) \cap Y\pr= \ov{v}_2^{-1}(U)  \hra U.
 \]
 Since $\phi$ is separated and \'etale, its fibre product with the composed immersion $Y_U\hra U\subset \bA^1_R$ is
 \[
  C\times_{\phi,\bA_{R}^1}Y_U=\wt{Y}_1 \sqcup \wt{Y}_2 \qq \text{ with } \qq  \wt{Y}_1\isoto Y_U.
 \]
   Let $W\subset \phi^{-1}(U)$ be an affine open such that $W\cap (C\times_{\phi,\bA_R^1}Y_U)=\wt{Y}_1$.
   Define $f\ce \phi|_W\colon W\ra U$, which has the desired properties by construction.
\epf

\section{ Torsors on a smooth affine relative curve}
\label{section-torsors on sm aff curves}
In this section, we establish a result regarding the triviality of torsors on a smooth affine relative curve. A similar result can also be found in \cite{Ces24}*{Theorem~4.4} and \cite{Fed23}*{Theorem~3}.
\bthmt[Section theorem]\label{triviality on sm rel. affine curves}
Let $R$ be a geometrically unibranch\footnote{We use the terminology of \SP{0BQ2}: each local ring is geometrically unibranch in the sense of \SP{0BPZ}; no Noetherian hypothesis is needed.}  
semilocal ring. Fix
\begin{itemize}
\item [-] a smooth affine $R$-curve $C$ and a section $s\in C(R)$, 
\item [-]an $R$-finite closed subscheme $Z\subset C$,
\item [-]a reductive $C$-group scheme $G$, 
\item [-]an $R$-algebra $A$, and a $G$-torsor $\cP$ over $C_A\ce C\times_RA$ that trivializes over $C_A\backslash Z_A$.
 \end{itemize}
 In each of the following cases
  \benumr
  \item \label{sec-thm-semilocal} $A$ is semilocal;
\item \label{sec-thm-totally isotropic} $s_A^*(G)$ is totally isotropic (cf. \S\ref{intro-Nisnevich purity});
  \eenum
  the pullback $s_A^*(\cP)$ is a trivial $s_A^*(G)$-torsor, where $s_A\in C_A(A)$ denotes the image of $s$. Moreover, {if $G$ descends to $R$}, then the {geometrically unibranch assumption} is redundant.
\ethmt

To prove Theorem \ref{triviality on sm rel. affine curves}, we first reduce to the case when $G$ descends to $R$; this is the only point where the geometrically unibranch assumption on $R$ is used. We then reduced to $C=\mathbb{A}^1_R$; see Lemmata \ref{equating red gps} and \ref{1st reduction of sec thm}. After these reductions, we use the geometry of the classifying stacks $\b{B}G$ and $\mathrm{Bun}G$.

\textbf{Reduction  to the case when $C=\mathbb{A}_R^1$ and $G$ is constant.}
We start with the following result for equating reductive group schemes; see also  \cite{Fed23}*{\S4.2}.

\blemt[\cite{Ces24}*{Lemma~4.3}]\label{equating red gps}
For a semilocal ring $B$, reductive $B$-group schemes $G_1$ and $G_2$ whose geometric $B$-fibres have the same type and whose radicals $\rad(G_1)$ and $\rad(G_2)$ are isotrivial,  maximal $B$-tori $T_1\subset G_1$ and $T_2\subset G_2$, and an ideal $I\subset B$, if there is an isomorphism of $B/I$-group schemes
\[
\text{$\iota\colon (G_1)_{B/I} \isoto (G_2)_{B/I}$  \q such that \q $\iota((T_1)_{B/I})=(T_2)_{B/I}$,}
\]
then there are a faithfully flat, finite, \'etale $B$-algebra $B\pr$, a section $s\colon B\pr\twoheadrightarrow B/I$, and an isomorphism of $B\pr$-group schemes $\iota\pr \colon (G_1)_{B\pr} \simeq (G_2)_{B\pr}$ such that $\iota\pr((T_1)_{B\pr})=(T_2)_{B\pr}$ and whose $s$-pullback is $\iota$.
\elemt

Now we state the promised reduction:

\blemt\label{1st reduction of sec thm}
    The proof of \Cref{triviality on sm rel. affine curves} reduces to the case when $C=\mathbb{A}_R^1$ and $G$ descends to $R$.
\elemt
\bpf
  Let $B$ be the semilocal ring of $C$ at the closed points of $\text{im}(s) \cup Z$; by \SP{0DQ1}, it is geometrically unibranch. By \cite{Guo24}*{Lemma~2.3}\footnote{To check the isotriviality, we may base change from $B$ to $B_{\mathrm{red}}\ce B/\sqrt{(0)}$, then the local rings of the latter are integral.}, the radicals of $G_B$ and $(s^*G)_B$ are isotrivial.
  By abuse of notation, we may view $s\colon B \twoheadrightarrow R$ as a section of the $R$-algebra $B$. As $B$ is semilocal, by \cite{SGA3II}*{exposé~XIV, corollaire~3.20}, $G_B$ admits a maximal $B$-torus $T_B$.
  Since the pullbacks of the pairs $(G_B,T_B)$ and $((s^*G)_B,(s^*T_B)_B)$ along $s$ are the same, by \Cref{equating red gps}, there are a faithfully flat, finite, \'etale $B$-algebra $B\pr$, a section $s\pr \colon B\pr \twoheadrightarrow R$ that lifts $s$, and a $B\pr$-isomorphism
  $$
  {\iota\colon (G_{B\pr},(T_B)_{B\pr}) \isoto \left((s^*G)_{B\pr},(s^*T_B)_{B\pr}\right)}
  $$
  whose $s$-pullback is the identity.
	  We spread out $\Spec B\pr \to \Spec B$ to obtain a finite \'etale cover $C\pr\to U $ of some affine open neighbourhood $U$ of $\text{im}(s) \cup Z$ in $C$.
  Shrinking $U$ if necessary, we may assume that the isomorphism $\iota$ is defined over $C\pr$. In both cases of \Cref{triviality on sm rel. affine curves} we may replace $C$ by $C\pr$, $Z$ by $C\pr\times_CZ$, $s$ by $s\pr$, and $\cP$ by $\cP|_{C\pr_A}$ to reduce to the case when $G$ is the base change of the reductive $R$-group scheme $s^*(G)$.

  Next, in order to apply gluing \cite{Ces22a}*{Lemma~7.1} to achieve that $C=\mathbb{A}_R^1$, we need to modify $C$ so that $Z$ embeds into $\mathbb{A}_R^1$.
		  For this, we first replace $Z$ by $Z\cup \text{im}(s)$ to assume that $s$ factors through $Z$. Then we apply Panin's `finite field tricks' \cite{Ces22a}*{Proposition~7.4} to obtain a finite morphism $\wt{C}\to C$ that is \'etale at the points in $\wt{Z}\ce \wt{C}\times_CZ$ such that $s$ lifts to $\wt{s}\in \wt{C}(R)$, and there is no finite field obstruction to embedding $\wt{Z}$ into $\mathbb{A}_R^1$ in the following sense: for every maximal ideal $\fm\subset R$,
\[
\# \bigl\{
z\in \wt{Z}_{k_{\fm}}\colon [k_z\colon k_{\fm}]=d \bigr\}< \# \bigl\{y\in \mathbb{A}_{k_{\fm}}^1\colon  [k_y\colon k_{\fm}]=d \bigr\} \q \text{for every} \q d\ge 1.
\]
Then, by \cite{Ces22a}*{Lemma~6.3}, there exist an affine open $C\prpr\subset \wt{C}$ containing $\text{im}(\wt{s})$ and a quasi-finite flat $R$-map $C\prpr\to \mathbb{A}_R^1$ that maps $\wt{Z}$ isomorphically onto a closed subscheme $Z\pr\subset \mathbb{A}_R^1$ with
    \[
	    \wt{Z}\simeq Z\pr \times_{\mathbb{A}_R^1}C\prpr.
    \]
(Actually, $C\prpr\to \mathbb{A}_R^1$ can be \'etale by shrinking $C\prpr$ around $\text{im}(\wt{s})$).
For both cases of \Cref{triviality on sm rel. affine curves}, since $\cP|_{C\prpr_A}$
is a $G$-torsor that trivializes over $C\prpr_A\backslash \wt{Z}_A$, we may use \cite{Ces22a}*{Lemma~7.1} to glue $\cP_{C\prpr_A}$ with the trivial $G$-torsor over $\mathbb{A}_A^1\backslash Z\pr_A$ to obtain a $G$-torsor $\cP\pr$ over $\mathbb{A}_A^1$ that trivializes over $\mathbb{A}_A^1 \backslash Z\pr_A$. Let $s\pr\in \mathbb{A}_R^1(R)$ be the image of $\wt{s}$; then ${s\pr}^{\ast}(\cP\pr)\simeq s^{\ast}(\cP)$.
It remains to replace $C$ by $\mathbb{A}_R^1$, $Z$ by $Z\pr$, $s$ by $s\pr$, and $\cP$ by $\cP\pr$.
\epf

\textbf{Torsors on $\bA^N$ and $\bP^1$.}
In an earlier version of this article, our analysis of torsors on $\mathbb{A}^1_R$ relied heavily on the geometry of affine Grassmannians developed by Fedorov (see \cite{Ces22}*{\S5.3} for a summary of and supplement to relevant techniques). 
In that approach, we made use of a slight variant of \cite{Ces22}*{Proposition~5.3.6} (which itself is a sharpening of \cite{Fed22b}*{Theorem~6}) to deduce the desired result. Later, \v{C}esnavi\v{c}ius and Fedorov  in \cite{CF23} developed an alternative approach by studying the \emph{open immersion} $\mathbf{B}G\hra \text{Bun}_G$, where $\mathbf{B}G$ denotes the classifying stack of $G$, and $\mathrm{Bun}_G$ denotes the moduli stack of $G$-torsors on $\bP^1$.
Since their results significantly simplify our argument, we have decided to cite them directly.

\bthm \label{Torsors-on-P1}
For a reductive group $G$ over a ring $A$ and a $G$-torsor $E$ on $\bP^1_A$,
\benumr
\item \label{semilocal-P1}$\mathrm{(}$\cite{CF23}*{Theorem~3.6}$\mathrm{)}$ if $A$ is semilocal, then $E|_{0_A}\simeq E|_{\infty_A}$;
\item \label{totallyisotropic-P1}$\mathrm{(}$\cite{CF23}*{Theorem~4.2}$\mathrm{)}$ if $G$ is totally isotropic and $E|_{\infty_A}$ is trivial, then $E|_{\bA^1_A}$ is trivial.
\eenum
\ethm

The above \Cref{Torsors-on-P1}~\ref{totallyisotropic-P1} immediately implies the following result, which was conjectured in \cite{Ces22}*{Conjecture~3.5.1} and settled later in \cite{Ces24}*{Corollary~2.2}.
\bthmt[\cite{Ces24}*{Cor.~2.2}, cf.~\cite{GL25}*{Thm.~3.5}]\label{triviality over relative affine space for totally isotropic}
  For a ring $R$ and a $\mathrm{totally}$ $\mathrm{isotropic}$ reductive $R$-group $G$, any $G$-torsor on $\mathbb{A}_R^N$ that is trivial away from some $R$-finite closed subscheme of $\mathbb{A}_R^N$ is trivial.
\ethmt
To deduce \Cref{triviality over relative affine space for totally isotropic} from \Cref{Torsors-on-P1} \ref{totallyisotropic-P1}, one first reduces it to the key case where $N=1$, then extends the torsor in question to a torsor over $\bP_R^1$ which trivializes along the infinity section.
Notice that the isotropicity assumption on $G$ is essential, see, e.g., \cite{Fed16} for counterexamples.

\bpp[Proof of \Cref{triviality on sm rel. affine curves}]
 Recall that $\cP$ is a $G$-torsor over $C_A$. 
 By \Cref{1st reduction of sec thm}, we may assume that $C=\mathbb{A}_R^1$ and $G$ is a reductive $R$-group scheme.
  After base change from $R$ to $A$, we may assume that $A=R$ (which need not be geometrically unibranch): in case \ref{sec-thm-semilocal}, $R$ is semilocal, and in case \ref{sec-thm-totally isotropic}, $G$ is totally isotropic.

\epp

For the totally isotropic case \ref{sec-thm-totally isotropic}, \Cref{triviality over relative affine space for totally isotropic} implies that $\cP$ is already trivial, and hence so is the pullback $s^*(\cP)$.

The semilocal case \ref{sec-thm-semilocal} follows directly from \cite{Ces24}*{Theorem~1.6}. 
Applied to $C=\bA_R^1$, $E=\cP$, and the section $s\in C(R)$, it proves the triviality of $s^*(\cP)$.

\section{ Torsors under a reductive group scheme over a smooth projective base}
\label{sect-torsor on sm proj base}

The main result of this section is the following:
\bthmt\label{torsors-Sm proj base}
   For a semilocal Pr\"{u}fer domain $R$, an $r\in R\backslash \{0\}$, an irreducible, smooth, projective $R$-scheme $X$, a finite subset $\mathbf{x} \subset X$ with semilocal ring $A\ce \sO_{X,\mathbf{x}}$, and a reductive $X$-group scheme $G$,
   \benumr
   \item \label{loc-gen-trivial-sm-proj} any generically trivial $G$-torsor over $A$ is trivial, that is,
   \[
   \ker\,(H^1_{\et}(A,G)\ra H^1_{\et}(\Frac A, G))=\{\ast\};
   \]
   \item \label{Nis-sm-proj}if $G_{A[\f{1}{r}]}$ is totally isotropic, then any generically trivial $G$-torsor over $A[\f{1}{r}]$ is trivial, that is,
   \[
   \tst \ker\,(H^1_{\et}(A[\f{1}{r}],G)\ra H^1_{\et}(\Frac A,G))=\{\ast\}.
   \]
  \eenum
\ethmt

  Part (i) proves the Grothendieck--Serre conjecture for reductive group schemes $G_A$ with a reductive model over a smooth projective compactification of $\Spec A$. Part (ii) gives the corresponding Nisnevich statement under the total-isotropy hypothesis. Specifically, if $R$ is a discrete valuation ring with uniformizer $r$ and $R\to A$ is a local homomorphism of local rings such that $r\in \mathfrak{m}_A\backslash \mathfrak{m}_A^2$, then case (ii) asserts that any generically trivial $G$-torsor over $A[\f{1}{r}]$ is trivial. Note that the isotropicity assumption on $G_A$ is essential (cf. \cite{Fed24}).

\bremt
An inspection of the proof below shows that \Cref{torsors-Sm proj base} still holds if $X$ is a flat projective $R$-scheme such that $X\backslash X^{\text{sm}}$ is $R$-fibrewise of codimension $\ge 2$ in $X$, $\mathbf{x}\subset X^{\text{sm}}$, and $G$ is a reductive $X^{\text{sm}}$-group scheme, where $X^{\text{sm}}$ denotes the smooth locus of $X\to \Spec R$.
\eremt

To prove \Cref{torsors-Sm proj base}, we first use \Cref{extend generically trivial torsors} and \Cref{Ces's Variant 3.7} to derive the following key result. This result reduces the proof of \Cref{torsors-Sm proj base} to studying torsors on a smooth affine relative curve.

\blemt\label{nicely spread out lem}
For a semilocal Pr\"{u}fer domain $R$ of finite Krull dimension, an irreducible, smooth, projective $R$-scheme $X$ of pure relative dimension $d> 0$, a finite subset $\mathbf{x} \subset X$, and a reductive $X$-group scheme $G$, the following assertions hold.
   \benumr
   \item \label{GS-nicely spread out} Given a generically trivial $G$-torsor $\cP$ over $A\ce \sO_{X,\mathbf{x}}$, there are
   \begin{itemize}
    \item  [-] a smooth, affine $A$-curve $C$, an $A$-finite closed subscheme $Z\subset C$, and a section $s\in C(A)$;
    \item [-] a reductive $C$-group scheme $\sG$ satisfying $s^{\ast}\sG\simeq G_A$ and a $\sG$-torsor $\cF$ such that $\cF|_{C\backslash Z}$ is trivial and $s^{\ast}\cF \simeq \cP$.
  \end{itemize}
  \item \label{Nis-nicely spread out} Given an $r\in R\backslash \{0\}$ and a generically trivial $G$-torsor $\wt{\cP}$ over $A[\f{1}{r}]$, there are
   \begin{itemize}
	    \item  [-] a smooth, affine $A$-curve $\wt{C}$, an $A$-finite closed subscheme $\wt{Z}\subset \wt{C}$, and a section $\wt{s}\in \wt{C}(A)$;
    \item [-] a reductive $\wt{C}$-group scheme $\wt{\sG}$ such that $\wt{s}^{\ast}\wt{\sG}\simeq G_A$, a $\wt{\sG}$-torsor $\wt{\cF}$ over $\wt{C}[\f{1}{r}]\ce \wt{C}\times_AA[\f{1}{r}]$ such that $\wt{\cF}|_{\wt{C}[\f{1}{r}]\backslash \wt{Z}[\f{1}{r}]}$ is trivial and $(\wt{s}|_{A[\f{1}{r}]})^{\ast}(\wt{\cF}) \simeq \wt{\cP}$.
  \end{itemize}
   \eenum
\elemt

\bpf
Since $X$ is projective over $R$, the finite set $\mathbf{x}$ is contained in an affine open. Indeed, for a projective embedding $X\hookrightarrow \bP_R^N$, homogeneous prime avoidance produces a positive-degree homogeneous form $f$ nonvanishing on $\mathbf{x}$, and $X\cap D_+(f)$ is affine. Thus $A=\sO_{X,\mathbf{x}}$ is defined.
Let $K$ be the fraction field of $R$.
By \Cref{extend generically trivial torsors}, $\cP$ (resp., $\wt{\cP}$) extends to a $G$-torsor $\cP_0$ (resp., $\wt{\cP_0}$) over an open neighbourhood $W\subset X$ of $\mathbf{x}$ (resp., an open neighbourhood $\wt{W}\subset X$ of $\Spec (A[\f{1}{r}])$) such that $(X\backslash W)_K$ (resp., $(X\backslash \wt{W})_K$) has codimension $\geq 3$ in $X_K$, and $(X\backslash W)_s$ (resp., $(X\backslash \wt{W})_s$) has codimension $\geq 2$ in $X_s$ for all $s\in \Spec (R)$.
Let $\textbf{z}\subset X$ be the set of maximal points of the $R$-fibres of $X$. The above codimension bounds imply that $\textbf{z}\subset W$ (resp., $\textbf{z}\subset \wt{W}$). By \Cref{geom}\ref{geo-iii}, the semilocal ring $\sO_{X,\textbf{z}}$, and hence also $\sO_{X,\textbf{z}}[\f{1}{r}]$, is a Pr\"{u}fer domain. By the Grothendieck--Serre result over semilocal Pr\"{u}fer domains (\Cref{G-S over semi-local prufer}), the generically trivial $G$-torsor $(\cP_0)|_{\sO_{X,\textbf{z}}}$ (resp., $(\wt{\cP_0})|_{\sO_{X,\textbf{z}}[\f{1}{r}]}$) is actually trivial. Thus there exists a closed subscheme $Y\subset X$ (resp., $\wt{Y}\subset X$) that avoids all the maximal points of the $R$-fibres of $X$ such that the restriction $(\cP_0)|_{X\backslash Y}$ (resp., $(\wt{\cP_0})|_{(X\backslash \wt{Y})[\f{1}{r}]}$) is trivial. Such a $Y$ (resp., $\wt{Y}$) is $R$-fibrewise of codimension $>0$ in $X$.
  Now, we treat the two cases \ref{GS-nicely spread out}--\ref{Nis-nicely spread out} separately.

For \ref{GS-nicely spread out}, by the above, $X\backslash W$ is $R$-fibrewise of codimension $\ge 2$ in $X$; \emph{a fortiori,} the same codimension bound holds for $Y\backslash W $ in $X$. Consequently, we can apply \Cref{Ces's Variant 3.7} to obtain an affine open $S\subset \mathbb{A}_{R}^{d-1}$, an affine open neighbourhood $ U\subset W$ of $\mathbf{x}$, and a smooth morphism $\pi\colon U\to S$ of pure relative dimension 1 such that $U\cap Y$ is $S$-finite.

Let $\tau\colon C\ce U\times_S\Spec A\to \Spec A$ be the base change of $\pi$ to $\Spec A$. Let $Z$ and $\cF$ be the pullbacks of $U\cap Y$ and $(\cP_0)|_U$ under $\mathrm{pr}_1\colon C\to U$, respectively. Then, via $\tau$, this $C$ is a smooth affine $A$-curve, $Z\subset C$ is an $A$-finite closed subscheme, and $\cF$ is a $\sG\ce \text{pr}_1^*(G_U)$-torsor that trivializes over $C\backslash Z$.
Finally, the diagonal in $C$ induces a section $s\in C(A)$ such that $s^{\ast}\sG\simeq G_A$ and $s^{\ast}\cF\simeq \cP$.

For \ref{Nis-nicely spread out},  since $\Spec A[\f{1}{r}] \subset \wt{W}$ and $\Spec (A[\f{1}{r}])$ consists of the points of $X[\f{1}{r}]\ce X\times_RR[\f{1}{r}]$ that specialize to some point of $\mathbf{x}$, no point of $(X\backslash \wt{W})[\f{1}{r}] =X[\f{1}{r}]\backslash \wt{W}[\f{1}{r}]$ specializes to any point of $\mathbf{x}$.
Therefore the closure of $\ov{(X\backslash \wt{W})[\f{1}{r}]}$ (in $X$) is disjoint from $\mathbf{x}$, and $\wt{W}\pr\ce X\backslash \ov{(X\backslash \wt{W})[\f{1}{r}]}$ is an open neighbourhood of $\mathbf{x}$.
    Notice that, since $\Spec R$ has finite Krull dimension, the topological space of $X$ is Noetherian.
    Since $(X\backslash \wt{W})[\f{1}{r}]$ is $R[\f{1}{r}]$-fibrewise of codimension $\ge 2$ in $X[\f{1}{r}]$, by \Cref{geom}\ref{geo-i} applied to the closures of the (finitely many) maximal points of $(X\backslash \wt{W})[\f{1}{r}]$, the closure $\overline{(X\backslash \wt{W})[\f{1}{r}]}=X \backslash \wt{W}\pr$ is $R$-fibrewise of codimension $\ge 2$ in $X$;
       \emph{a fortiori}, the same holds for $\wt{Y}\backslash \wt{W}\pr$ in $X$. Consequently, we can apply \Cref{Ces's Variant 3.7} to obtain an affine open $\wt{S}\subset \mathbb{A}_{R}^{d-1}$, an affine open neighbourhood $ \wt{U}\subset \wt{W}\pr$ of $\mathbf{x}$, and a smooth morphism $\wt{\pi}\colon \wt{U}\to \wt{S}$ of pure relative dimension 1 such that $\wt{U}\cap \wt{Y}$ is $\wt{S}$-finite. Since $\wt{U}[\f{1}{r}]\subset \wt{W}\pr[\f{1}{r}] =\wt{W}[\f{1}{r}]$, we can form the restriction $(\wt{\cP_0})|_{\wt{U}[\f{1}{r}]}$.

       Let $\wt{\tau}\colon \wt{C}\ce \wt{U}\times_{\wt{S}}\Spec A\to \Spec A$ be the base change of $\wt{\pi}$ to $\Spec A$. Let $\wt{Z}$ be the pullback of $\wt{U}\cap \wt{Y}$ under $\mathrm{pr}_1\colon \wt{C}\to \wt{U}$. Let $\wt{\cF}$ be the pullback of $(\wt{\cP_0})|_{\wt{U}[\f{1}{r}]}$ under $\mathrm{pr}_1\colon \wt{C}[\f{1}{r}]\to \wt{U}[\f{1}{r}]$. Then, via $\wt{\tau}$, this $\wt{C}$ is a smooth affine $A$-curve, $\wt{Z}\subset \wt{C}$ is an $A$-finite closed subscheme, and $\wt{\cF}$ is a $\wt{\sG}\ce \text{pr}_1^{\ast}(G_{\wt{U}})$-torsor over $\wt{C}[\f{1}{r}]$ that trivializes over $\wt{C}[\f{1}{r}]\backslash \wt{Z}[\f{1}{r}]$. Finally, the diagonal in $\wt{C}$ induces a section $\wt{s}\in \wt{C}(A)$ such that $\wt{s}^{\ast}\wt{\sG}\simeq G_A$ and $(\wt{s}|_{A[\f{1}{r}]})^{\ast}(\wt{\cF})\simeq \wt{\cP}$.\qedhere
\epf

\bpp[Proof of \Cref{torsors-Sm proj base}]
A standard argument involving \Cref{geom}\ref{limitarg} reduces to the case when $R$ has finite Krull dimension.
Now, let $\cP$ (resp., $\wt{\cP}$) be a generically trivial $G$-torsor over $A\ce \sO_{X,\mathbf{x}}$ (resp., over $A[\f{1}{r}]$).
 Let $d$ be the relative dimension of $X$ over $R$.
 If $d=0$, then $A$ and $A[\f{1}{r}]$ are semilocal Pr\"{u}fer domains, so, by the Grothendieck--Serre result over semilocal Pr\"{u}fer domains (\Cref{G-S over semi-local prufer}), the torsors $\cP$ and $\wt{\cP}$ are trivial. Hence we may assume that $d>0$. Then, by \Cref{nicely spread out lem}, there are a smooth, affine $A$-curve $C$ (resp., $\wt{C}$), an $A$-finite closed subscheme $Z\subset C$ (resp., $\wt{Z}\subset \wt{C}$), a section $s\in C(A)$ (resp., $\wt{s}\in \wt{C}(A)$), a reductive $C$-group scheme $\sG$ (resp., a reductive $\wt{C}$-group scheme $\wt{\sG}$) with $s^{\ast}\sG\simeq G_A$ (resp., $\wt{s}^{\ast}\wt{\sG}\simeq G_A$),
 \begin{itemize}
   \item [-]  a $\sG$-torsor $\cF$ over $C$ that trivializes over $C\backslash Z$ such that $s^{\ast}\cF \simeq \cP$,
   \item [-] resp., a $\wt{\sG}$-torsor $\wt{\cF}$ over $\wt{C}[\f{1}{r}]$ that trivializes over $\wt{C}[\f{1}{r}]\backslash \wt{Z}[\f{1}{r}]$ such that $(\wt{s}|_{A[\f{1}{r}]})^{\ast}(\wt{\cF}) \simeq \wt{\cP}$.
 \end{itemize}
 By \Cref{triviality on sm rel. affine curves} \ref{sec-thm-semilocal}, the $G$-torsor $s^{\ast}\cF \simeq \cP$ is trivial. By \Cref{triviality on sm rel. affine curves} \ref{sec-thm-totally isotropic}, in case $(\wt{s}|_{A[\f{1}{r}]})^{\ast}(\wt{\sG})\simeq G_{A[\f{1}{r}]}$ is totally isotropic, the $G_{A[\f{1}{r}]}$-torsor $(\wt{s}|_{A[\f{1}{r}]})^{\ast}(\wt{\cF}) \simeq \wt{\cP}$ is trivial.
 \QED
\epp

\section{ Torsors under a constant reductive group scheme}
\label{sect-torsor under constant redu}
In this section, we establish the first main result of this paper, namely, the Grothendieck--Serre conjecture and a version of the Nisnevich conjecture for `constant' reductive group schemes. Our proof relies on a variant of Lindel's Lemma (\Cref{variant of Lindel's lem}) and gluing techniques to reduce to the case proved in \Cref{torsors-Sm proj base}.
\bthmt\label{G-S for constant reductive gps}
   For a semilocal Pr\"ufer domain $R$, an $r\in R\backslash \{0\}$, an irreducible affine $R$-smooth scheme $X$, a finite subset $\b{x}\subset X$, and a reductive $R$-group scheme $G$,
   \benumr
   \item\label{G-S for constant reductive gps i} any generically trivial $G$-torsor over $A\ce \sO_{X,\mathbf{x}}$ is trivial, that is,
    \[
\mathrm{ker}\left(H^1_{\et}(A,G)\to H^1_{\et}(\Frac A,G)\right)=\{*\};
    \]
   \item\label{G-S for constant reductive gps ii} if $(G_{R[\f{1}{r}]})^{\mathrm{ad}}$ is totally isotropic, then any generically trivial $G$-torsor over $A[\f{1}{r}]$ is trivial, that is,
    \[
    \tst \mathrm{ker}\left(H^1_{\et}(A[\f{1}{r}],G)\to H^1_{\et}(\Frac A,G)\right)=\{*\}.
    \]
  \eenum
\ethmt


\bpf
Let $\cP$ (resp., $\wt{\cP}$) be a generically trivial $G$-torsor over $A$ (resp., over $A[\f{1}{r}]$). We can assume that $\cP$ is defined on the whole $X$ (resp. $\wt{\cP}$ is defined on the whole $X[\f{1}{r}]\ce X\times_RR[\f{1}{r}]$) by shrinking $X$ around $\mathbf{x}$. Let $d$ be the relative dimension of $X$ over $R$.    Let $i_0\colon X\to X\times_R \bA_R^N$ denote the zero section. As pointed out by $\check{\mathrm{C}}$esnavi$\check{\mathrm{c}}$ius, for a large $N$, we may replace $X$ by $X\times_R \bA_R^N$, $\cP$ by $\cP_{X\times_R \bA_R^N}$, and $\mathbf{x}$ by $i_0(\mathbf{x})$ so that $d>\#\, \mathbf{x}$: in fact, the triviality of $\cP_{X\times_R \bA_R^N}$ (resp., $\wt{\cP}_{X[\f{1}{r}]\times_R \bA_R^N}$) around $i_0(\mathbf{x})$ pulls back along $i_0$ to that of $\cP$ (resp., $\wt{\cP}$) around $\mathbf{x}$.

 Our goal is to prove that $\cP|_A$ (resp., $\wt{\cP}|_{A[\f{1}{r}]}$) is trivial. 
 By a standard limit argument \Cref{geom}\ref{limitarg}, we are reduced to the case when $R$ has finite Krull dimension. Using specialization, we can further assume that each point of $\mathbf{x}$ is closed in its corresponding $R$-fibre of $X$, although it does not necessarily lie in a closed $R$-fibre of $X$.
 If $d=0$, then $A$ (resp., $A[\f{1}{r}]$) is a semilocal Pr\"{u}fer domain. By \Cref{G-S over semi-local prufer}, the torsor $\cP|_A$ (resp., $\wt{\cP}|_{A[\f{1}{r}]}$) is then trivial. Thus we may assume that $d>0$ for what follows. Denote by $\textbf{y}$ the set of maximal points of the $R$-fibres of $X$.
\bclt
No point of $\b{x}$ specializes to any point of $\b{y}$, that is, $\overline{\mathbf{x}} \cap \b{y} =\emptyset$.
\eclt
  \bpf [Proof of the claim] Let $\pi\colon X\to S\ce \Spec R$ be the structural morphism.
  Assuming the claim is false, let $x\in \mathbf{x}$ specialize to $y\in \textbf{y}$. Then, by \Cref{geom}\ref{geo-i}, we have
  \[
 \dim \ov{\{x\}}_{\pi(x)} = \dim \ov{\{x\}}_{\pi(y)},
  \]
  which is $\ge \dim \overline{\{y\}}_{\pi(y)}=d$ (because $y$ is a maximal point in the fibre $ \pi^{-1}(\pi(y))$ which has pure dimension $d$). Since $\dim \pi^{-1}(\pi(x))=d>0$, the point $x$ cannot be a closed point of the fibre $ \pi^{-1}(\pi(x))$, a contradiction.
  \epf

   Using \Cref{geom}\ref{geo-iii} again, we see that the semilocal ring $\sO_{X,\textbf{y}}$ (and hence also $\sO_{X,\textbf{y}}[\frac{1}{r}]$) is a Pr\"{u}fer domain. Therefore, by \Cref{G-S over semi-local prufer}, the $G$-torsor $\cP|_{\sO_{X,\textbf{y}}}$ (resp., $\wt{\cP}|_{\sO_{X,\textbf{y}}[\frac{1}{r}]}$) is trivial. Using the above claim and prime avoidance \SP{00DS}, we can find an element $a\in \Gamma(X,\sO_X)$ such that $Y\ce V(a)\subseteq X$ contains $\mathbf{x}$, avoids $\textbf{y}$, and the restriction $\cP|_{X\setminus Y}$ (resp., $\wt{\cP}|_{(X\setminus Y)[\frac{1}{r}]}$) is trivial.

Since $d>\#\, \mathbf{x}$, we can apply \Cref{variant of Lindel's lem} to obtain an affine open neighbourhood $W\subset X$ of $\mathbf{x}$, an affine open subscheme $ U\subset \mathbb{A}_R^d$, and an \'etale $R$-morphism $f\colon W\to U$ such that the restriction $f|_{W\cap Y}$ is a closed immersion and $f$ induces a Cartesian square
\begin{equation*}
\begin{tikzcd}
W\cap Y \arrow[r, hook]  \arrow[d, equal]
 & W \arrow[d, "f"] \\
W\cap Y \arrow[r, hook]
& U.
\end{tikzcd}
\end{equation*}
 Applying $(-)\times_RR[\f{1}{r}]$ yields a similar Cartesian square. Using gluing \cite{Ces22a}*{Lemma~7.1}, we can treat the two cases as follows.
\benumr
\item We may (non-canonically) glue $\cP|_{W}$ with the trivial $G$-torsor over $U\backslash f(W\cap Y)$ to descend $\cP|_{W}$ to a $G$-torsor $\cQ$ over $U$ that trivializes over $U\backslash f(W\cap Y)$. Since $U$ has a smooth, projective compactification $\mathbb{P}_R^d$, we may apply \Cref{torsors-Sm proj base}~\ref{loc-gen-trivial-sm-proj} to deduce that $\cQ|_{\sO_{U,f(\mathbf{x})}}$ is trivial. Therefore, $\cP|_A=\cP|_{\sO_{W,\mathbf{x}}}$
    is trivial, as desired.
\item We may (non-canonically) glue $\wt{\cP}|_{W[\f{1}{r}]}$ with the trivial $G$-torsor over $(U\backslash f(W\cap Y))[\f{1}{r}]$ to descend $\wt{\cP}|_{W[\f{1}{r}]}$ to a $G$-torsor $\wt{\cQ}$ over $U[\f{1}{r}]$ that trivializes over $U[\f{1}{r}]\backslash f(W\cap Y)[\f{1}{r}]$. Since $U$ has a smooth, projective compactification $\mathbb{P}_R^d$, we may apply \Cref{torsors-Sm proj base}~\ref{Nis-sm-proj} to conclude that $\wt{\cQ}|_{\sO_{U,f(\mathbf{x})}[\f{1}{r}]}$ is trivial. Therefore, $\wt{\cP}|_{A[\f{1}{r}]}=\wt{\cP}|_{\sO_{W,\mathbf{x}}[\f{1}{r}]}$
    is trivial, as desired. \qedhere
\eenum
\epf

\begt
Let $D$ be an Azumaya algebra over a discrete valuation ring $R$ (or, more generally, a valuation ring $R$), and let $A$ be an ind-smooth, local $R$-algebra with fraction field $K$. Then, a unit of $A$ is a reduced norm from $D_A\ce D\otimes_RA$ if and only if it is a reduced norm from $D_K\ce D\otimes_RK$. When $R$ is an infinite field, this result was obtained in \cite{CTO92}*{corollaire 5.2} via a reduction trick from $\SL(D)$ to $\SL_2(D)$ (which is totally isotropic) and an application of the Grothendieck--Serre conjecture for the latter. We now justify the assertion in our setting.

\Cref{G-S for constant reductive gps}~\ref{G-S for constant reductive gps i} gives a direct proof.. Consider the following short exact sequence of reductive $R$-groups:
\[
1\to \SL(D) \to \GL(D) \to \mathbb{G}_{m,R}\to 1.
\]
For any local $R$-algebra $A$, we have $H^1_{\et}(A,\GL(D))=\{*\}$: indeed, a $\GL(D)$-torsor is a right $D_A$-module $P$ that is \'etale-locally free of rank one; choosing a generator of $P/\mathfrak{m}_AP$ over $D_A/\mathfrak{m}_AD_A$ and lifting it yields an isomorphism $D_A\simeq P$ by Nakayamas lemma. From the associated long exact cohomology sequence, we obtain the bijections 
\[
A^{\times}/\mathrm{Nrd}(D_A^{\times}) \simeq H^1_{\et}(A,\SL(D)), \qq \text{and} \qq K^{\times}/\mathrm{Nrd}(D_K^{\times}) \simeq H^1_{\et}(K,\SL(D)).
\]
 Thus, the injectivity of the homomorphism
\[
A^{\times}/\mathrm{Nrd}(D_A^{\times})\to K^{\times}/\mathrm{Nrd}(D_K^{\times})
\]
follows directly from \Cref{G-S for constant reductive gps}~\ref{G-S for constant reductive gps i}, applied to the reductive $R$-group $G=\SL(D)$.
\eegt

As a corollary of \Cref{G-S for constant reductive gps}~\ref{G-S for constant reductive gps i}, we obtain the following non-Noetherian generalisation of \cite{Ces22a}*{Corollary 9.6} (which in turn generalises the field case obtained in \cite{CTO92}*{corollaire 2.8}).
\bcort
 For a semilocal ring $A$ that is ind-smooth over a Pr\"ufer ring with $2 \in A^\times$, we have
\[
H^1(A, \mathrm{SO}_n) \hookrightarrow H^1(\mathrm{Frac}(A), \mathrm{SO}_n) 
\quad \text{and} \quad 
H^1(A, \mathrm{O}_n) \hookrightarrow H^1(\mathrm{Frac}(A), \mathrm{O}_n)
\]
for all $n \geq 1$. Moreover, two nondegenerate quadratic forms over $A$ are isomorphic provided that they become isomorphic over $\mathrm{Frac}(A)$.
\ecort
\bpf
The proof is parallel to that of \cite{Ces22a}*{Corollary~9.6}: one only replaces the Grothendieck--Serre input used there by \Cref{G-S for constant reductive gps}\ref{G-S for constant reductive gps i}.
\epf

\begin{appendix}

\section{Grothendieck--Serre on a semilocal Pr\"{u}fer domain}
\label{section-G-S on semilocal prufer}

We generalise the Grothendieck--Serre theorems in \cite{Guo22} and \cite{Guo24} to semilocal Pr\"{u}fer rings.
This theorem is a key input for our main results, Theorems \ref{torsors-Sm proj base} and \ref{G-S for constant reductive gps}, and is stated as follows.
\bthm \label{G-S over semi-local prufer}
For a semilocal Pr\"{u}fer domain $R$ and a reductive $R$-group scheme $G$, we have
\[
   \text{$\mathrm{ker}\left(\text{H}_{\et}^1(R,G)\to \text{H}_{\et}^1(\Frac R,G)\right)=\{*\}$}.
\]
\ethm

Using a standard limit argument involving \Cref{geom}\ref{limitarg}, we reduce to the case where $R$ has finite Krull dimension.
We then argue by induction on $n=\dim R$.
Choose $\varpi\in R$ with $V(\varpi)=\operatorname{Max}R$.
The induction hypothesis applies over $R[\frac{1}{\varpi}]$, while \cite{Guo24}*{Theorem~1.3} applies over the $\varpi$-adic completions of the local rings at the maximal ideals.
The remaining step is patching.
\bpp[Setup]\label{pd-setup}
Let $R$ be a semilocal Pr\"ufer domain of finite Krull dimension. 
Let $(\fm_i)_{i=1}^r$ be all the maximal ideals of $R$ with the local rings $R_i\ce R_{\fm_i}$.
Let $\varpi\in R$ be such that $V(\varpi)=\{\fm_i\}_{i=1}^r$.
Let $\wh{R}$ (resp., $\wh{R_i}$) denote the $\varpi$-adic completion of $R$ (resp., of $R_i$).
Then $\wh{R_i}$ is an $\varpi$-adically complete valuation ring of rank $1$, and we have an isomorphism of topological rings $\wh{R}\simeq \prod_{i=1}^r\wh{R_i}$. 
Set $\wh{K}_i\ce \Frac \wh{R_i}=\wh{R_i}[\f{1}{\varpi}]$; it is a $\varpi$-adically complete valued field.
Topologize $R[\f{1}{\varpi}]$ by declaring $\{\mathrm{im}(\varpi^nR\ra R[\f{1}{\varpi}])\}_{n\ge 1}$ to be a fundamental system of open neighbourhoods of $0$.
This is the topology induced from the $\varpi$-adic topology on $R$, so \cite{BC22}*{Lemma~2.1.11} identifies the completion
\[
\tst \text{$R[\f{1}{\varpi}]\to \wh{R}[\f{1}{\varpi}]\simeq  \prod_{i=1}^r \wh{R_i}[\f{1}{\varpi}]=\prod_{i=1}^r \wh{K}_i$.}
\]
We write $R_\varpi\ce R[\frac{1}{\varpi}]$, $F\ce\prod_i\wh{K}_i$, and $\mathcal O\ce\prod_i\wh{R_i}$.

In the sequel, let $G$ be a reductive $R$-group scheme.

\epp
We use the Beauville--Laszlo patching form of \cite{Guo24}*{Proposition~4.5}.
\blem \label{double cosets}
      The $G$-torsors on $R$ that trivialize both on $R[\f{1}{\varpi}]$ and on $\prod_{i=1}^r\wh{R_i}$ are in bijection with the double cosets 
      \[
      \tst \x{$\mathrm{im}\bigl(G(R[\f{1}{\varpi}])\ra \prod_{i=1}^r G(\wh{K}_i)\bigr)\backslash \prod_{i=1}^r G(\wh{K}_i)/ \prod_{i=1}^rG(\wh{R_i})$.}
      \]
      \elem
\bpf

Applying \cite{BC22}*{Lemma~2.2.11(b)} with $A=R$, $A\pr=\wh{R}$, and $t=\varpi$, we obtain the Beauville--Laszlo patching equivalence for $G$-torsors.
\par
The hypotheses are satisfied because $R/\varpi^mR\simeq\wh{R}/\varpi^m\wh{R}$ for $m>0$, while $\varpi$ is a nonzerodivisor on $R$ and on $\wh{R}\simeq\prod_i\wh{R_i}$.
Since $G$ is affine and flat over $R$, the resulting patching equivalence means that a torsor trivial on $R[\f{1}{\varpi}]$ and on $\wh{R}$ is obtained from a gluing element of
\[
\tst {G(\wh{R}[\f{1}{\varpi}])=\prod_iG(\wh{K}_i).}
\]
Changing the two trivializations translates this element by left multiplication via the image of $G(R[\f{1}{\varpi}])$ and by right multiplication via $\prod_iG(\wh{R_i})$.
Quotienting  by these changes is exactly a double coset as displayed in the statement.

\epf
If a $G$-torsor over $R$ trivializes over $R[\f{1}{\varpi}]$, then it also trivializes over $\prod_{i=1}^r\wh{R_i}$ by \cite{Guo24}*{Theorem~1.3}.
Consequently, \Cref{double cosets} reduces the patching step to the following product formula.
 \bprop\label{decomp-gp}
      Under Setup~\ref{pd-setup}, we have
      \[
      \tst \prod_{i=1}^r G(\wh{K}_i)=\mathrm{im}\bigl(G(R[\f{1}{\varpi}])\ra \prod_{i=1}^r G(\wh{K}_i)\bigr)\cdot \prod_{i=1}^rG(\wh{R_i}).
      \]
\eprop
The proof uses the following approximation lemmata.
\blem\label{lift-tor}
    Let $S$ be a semilocal scheme and let $G$ be a reductive $S$-group scheme.
    Let $\mtg$ be the scheme of maximal tori of $G$.
   If $\# k_s\geq \dim (G_{k_s}/Z_{k_s})$ for every closed point $s\in S$, where $Z_{k_s}\subset G_{k_s}$ is the centre, then the following natural map is surjective
      \[
        \tst  \mtg(S)\surjects \prod_{s\in S\mathrm{\, closed}}\mtg(k_s).
      \]
      \elem
\bpf
By \cite{SGA3II}*{exposé~XII, théorème~4.7~c)}, we may replace $G$ by its adjoint quotient.
Write $S=\Spec A$, and let $(T_s)_s$ be a tuple of maximal tori over the closed fibres.
By \cite{SGA3II}*{exposé~XIV, théorèmes~3.9 and~3.18}, it suffices to lift the Cartan subalgebras $\Lie(T_s)\subset\Lie(G_s)$.
By the cardinality assumption and \cite{Bar67}*{Theorem~1}, for each closed $s$ choose $a_s\in\Lie(T_s)$ such that $\Lie(T_s)=\bigcup_n\ker(\mathrm{ad}(a_s)^n)$.
Then $a_s$ is regular in $\Lie(G_s)$ by \cite{SGA3II}*{exposé~XIII, corollaire~5.7}.
The finite local freeness of $\Lie(G)$ and the surjection $A\to\prod_s k_s$ allow us to lift $(a_s)_s$ to a section $a\in\Gamma(S,\Lie(G))$.
The regular locus is open (\emph{op. cit.}) and contains all closed points of $S$, hence $a$ is a regular section.
By \cite{SGA3II}*{exposé~XIII, corollaire~5.7 and exposé~XIV, théorème~3.18}, $\mathfrak c\ce\mathrm{Nil}(a)=\bigcup_{n\ge 1}\ker(\mathrm{ad}(a)^n)$ is the Cartan subalgebra attached to $a$, with closed fibres $\Lie(T_s)$.
The maximal torus attached to $\mathfrak c$ lifts the tuple $(T_s)_s$.
\epf

\blem\label{approx-tori-and-normal-subgroup}
Under Setup~\ref{pd-setup}, the following hold.
\benumr
\item The image of $\mtg(R_\varpi)\to\mtg(F)$ is dense in the $\varpi$-adic topology.
\item For every maximal $\mathcal O$-torus $T_{\mathcal O}\subset G_{\mathcal O}$ and every open neighbourhood $W\subset G(F)$ of $1$, there exist a maximal $R$-torus $T_0\subset G$ and a $g\in W$ such that $(T_0)_F=g(T_{\mathcal O})_Fg^{-1}$.
\item If $C\subset G(F)$ is the closure of the image of $G(R_\varpi)$ in $G(F)$, then $C$ contains an open normal subgroup of $G(F)$.
\eenum
\elem
\bpf

\noindent\textup{(i).}
For $m\ge 0$, let $\mathrm{Cauchy}^{\ge m}(R_\varpi)$ be the ring of Cauchy sequences $(a_N)_{N\ge m}$ in $R_\varpi$ for the topology of \S\ref{pd-setup}, with termwise addition and multiplication.
The transition map from level $m$ to level $m+1$ forgets the first term.
Denote
\[
   {\tst A_{\infty}\ce
   \varinjlim_{m\geq 0}\mathrm{Cauchy}^{\geq m}(R_\varpi)}
\]
so an element of $A_\infty$ is a Cauchy tail modulo finitely many initial terms.
The limit map $A_{\infty}\to F$ is surjective.
Let $J\ce\ker(A_\infty\to F)$.
If $b\in J$, then $1+b$ is represented by a tail of units in $R_\varpi$, so $1+b$ is a unit in $A_\infty$.
Thus $J\subset \Jac(A_\infty)$, and the maximal ideals of $A_\infty$ are the kernels of the projections $A_\infty\to\wh{K}_i$.
In particular, $A_\infty$ is semilocal with residue fields $\wh{K}_i$.

By \cite{SGA3II}*{exposé~XIV, théorème~6.1}, the scheme $\mtg$ is affine and finitely presented, so
\[
   \tst
   \mtg(A_{\infty})=
   \varinjlim_{m\geq 0}\mtg\bigl(\mathrm{Cauchy}^{\geq m}(R_\varpi)\bigr).
\]
The fields $\wh{K}_i$ are infinite.
By \Cref{lift-tor}, every point of $\mtg(F)=\prod_i\mtg(\wh{K}_i)$ lifts to $\mtg(A_\infty)$.
Such a lift is represented by a Cauchy tail of points of $\mtg(R_\varpi)$.

\noindent\textup{(ii).}
Let $T\ce (T_{\mathcal O})_F$.
Consider the conjugacy morphism
\[
   \tst
   G_F\longrightarrow \mtg_F,\qquad g\longmapsto gTg^{-1}.
\]
Here $G_F/N_{G_F}(T)=\mtg_F$.
By \cite{Ces15d}*{4.3(a) and 2.8(2)}, the image of $W$ in $\mtg(F)$ contains an open neighbourhood of $T$.
Moreover, $\mtg(\mathcal O)\subset\mtg(F)$ is open by \cite{Guo24}*{Lemma~3.5(iii)}.
By the density assertion, choose a $T_\varpi\in\mtg(R_\varpi)$ whose image in $\mtg(F)$ lies in the intersection.
Choose $T_{\mathcal O}\pr\in\mtg(\mathcal O)$ with the same image in $\mtg(F)$.
By \cite{BC22}*{Proposition~2.2.12}
\[
   \tst
   \mtg(R)\simeq \mtg(R_\varpi)\times_{\mtg(F)}\mtg(\mathcal O),
\]
which glues $T_\varpi$ and $T_{\mathcal O}\pr$ to a maximal torus $T_0\subset G$ with $(T_0)_F=gTg^{-1}$ for some $g\in W$.

\noindent\textup{(iii).}
Let $C\subset G(F)$ be the closure of the image of $G(R_\varpi)$.
We use \cite{Guo24}*{Lemma~3.17 and Proposition~3.19} in product form over $F$.
Let $T\pr$ be an $R_\varpi$-torus, and let $B/R_\varpi$ be a finite Galois splitting algebra of $T\pr$.
Set $B_\infty\ce B\otimes_{R_\varpi}A_\infty$ and $B_F\ce B\otimes_{R_\varpi}F$.
For $J=\ker(A_\infty\to F)$, finite flatness of $B$ implies $\ker(B_\infty\to B_F)=B\otimes_{R_\varpi}J\subset\Jac(B_\infty)$.
Hence $B_\infty^\times\to B_F^\times$ is surjective.
Since $T\pr_B$ is split, $T\pr(B_\infty)\to T\pr(B_F)$ is surjective.
As $T\pr$ is affine and finitely presented, each lift in $T\pr(B_\infty)$ is represented by a Cauchy tail in $T\pr(B)$.
Applying $\mathrm{Nm}_{B/R_\varpi}$ termwise and passing to the limit, we obtain
\[
\mathrm{Nm}_{B_F/F}(T\pr(B_F))\subset \overline{T\pr(R_\varpi)}.
\]

Fix a maximal $F$-torus $T\subset G_F$.
Let $L/F$ be a finite Galois splitting algebra of $T$ and $U_T\subset T(F)$ the image of $N_{L/F}\colon T(L)\to T(F)$, which is open by \cite{Guo24}*{Lemma~3.16(ii)}.
We first show that $U_T\subset C$.
Let $u\in U_T$ and let $\Omega$ be an open neighbourhood of $u$ in $G(F)$.
Choose an open neighbourhood $W$ of $1$ in $G(F)$ such that $huh^{-1}\in\Omega$ for all $h\in W$.
From the two assertions proved above, take $h\in W$ and a maximal $R_\varpi$-torus $T_\varpi\subset G_{R_\varpi}$ with $(T_\varpi)_F=hTh^{-1}$.
By the norm inclusion above, $hU_Th^{-1}\subset \overline{T_\varpi(R_\varpi)}\subset C$.
Every neighbourhood $\Omega$ of $u$ therefore meets $C$.
Since $C$ is closed, $u\in C$, and hence $U_T\subset C$.
The choice of $T$ was arbitrary, so this holds for every maximal $F$-torus.

Let $T^{\mathrm{reg}}\subset T$ be the regular open locus.
The norm morphism $\mathrm{Res}_{L/F}(T_L)\to T$ is dominant.
Since $L$ splits $T$ and $F$ is a product of infinite fields, $T(L)=\mathrm{Res}_{L/F}(T_L)(F)$ is Zariski dense in $\mathrm{Res}_{L/F}(T_L)$.
Thus $U_T\cap T^{\mathrm{reg}}(F)\ne\emptyset$.
Choose $t\in U_T\cap T^{\mathrm{reg}}(F)$.
By \cite{SGA3II}*{exposé~XIII, proposition~2.2}, the map
\[
   \tst
   f\colon G_F\times T\longrightarrow G_F,\qquad (g,t\pr)\longmapsto gt\pr g^{-1}
\]
is smooth at $(1,t)$.
By the criterion for openness  \cite{Guo24}*{Lemma~3.6(i)} and the smoothness of $f$ at the $G(F)$-translates of $(1,t)$, there is an open neighbourhood $V\subset U_T\cap T^{\mathrm{reg}}(F)$ of $t$ such that
\[
   \tst
   E\ce \bigcup_{g\in G(F)}gVg^{-1}
\]
is open in $G(F)$.
For each $g\in G(F)$, the same inclusion for $gTg^{-1}$ shows that $gU_Tg^{-1}\subset C$.
Thus $E\subset C$.
Let $N$ be the subgroup generated by $E$.
Since $E$ is open and stable under conjugation, $N$ is open and normal in $G(F)$.
Since $C$ is a subgroup containing $E$, it also contains $N$.
\epf

\blem[\cite{Guo24}*{Proposition~4.5}]\label{closed-decomp-gp}
Under Setup~\ref{pd-setup}, let $C\subset G(F)$ be the closure of the image of $G(R_\varpi)$.
Then $G(F)=C\cdot G(\mathcal O)$.
\elem
\bpf
For an $R$-torus $T$, \Cref{double cosets} identifies $\mathrm{im}\bigl(T(R_\varpi)\to T(F)\bigr)\backslash T(F)/T(\mathcal O)$ with the set of $T$-torsors on $R$ that are trivial on $R_\varpi$ and on $\mathcal O$.
These torsors are generically trivial, and the map $H^1_{\et}(R,T)\to H^1_{\et}(\Frac R,T)$ is injective by \Cref{G-S type results for mult type}.
Consequently,
\begin{equation}\label{torus-factorization}
   \tst T(F)=
   \mathrm{im}\bigl(T(R_\varpi)\to T(F)\bigr)\cdot T(\mathcal O).
\end{equation}

Let $N\subset C$ be the open normal subgroup supplied by \Cref{approx-tori-and-normal-subgroup}~\textup{(iii)}.
Suppose first that $G_F$ has no proper parabolic.
Let $Z\ce\rad(G)$ and $H\ce G/Z$.
Over the valuation rings $\wh{R_i}$, we have $H(F)=H(\mathcal O)$ by \cite{Guo24}*{Proposition~4.3(e)}.
The injectivity of the map $H^1(\mathcal O,Z)\to H^1(F,Z)$ follows from \Cref{G-S type results for mult type} over the valuation rings $\wh{R_i}$.
The exact cohomology sequence attached to $1\to Z\to G\to H\to 1$ therefore reduces an element of $G(F)$, modulo $G(\mathcal O)$, to an element of $Z(F)$.
By \eqref{torus-factorization} for the $R$-torus $Z$, we have
\[
   \tst Z(F)=
   \mathrm{im}\bigl(Z(R_\varpi)\to Z(F)\bigr)\cdot Z(\mathcal O)
   \subset C\cdot G(\mathcal O).
\]

It follows that $G(F)\subset C\cdot G(\mathcal O)$.

Now suppose that $G_F$ has a proper parabolic.
By \cite{SGA3IIInew}*{exposé~XXVI, corollaire~3.5}, the scheme of parabolics is proper over $\mathcal O$.
Thus a minimal parabolic of $G_F$ extends to a parabolic $P\subset G_{\mathcal O}$.
Choose an opposite parabolic $Q$, a common Levi $L$, and write $U$ and $U_Q$ for the two unipotent radicals.
Let $S\subset L$ be a maximal split central torus.
The inclusion $N\subset C$ makes $C\cap G(\mathcal O)$ an open neighbourhood of $1$ in $G(F)$.
Choose a maximal $\mathcal O$-torus $T\subset L$ containing $S$, obtained from a maximal torus of $L/S$.
By the maximal-torus approximation in \Cref{approx-tori-and-normal-subgroup} with neighbourhood $C\cap G(\mathcal O)$, choose $g\in C\cap G(\mathcal O)$ and a maximal $R$-torus $T_0\subset G$ with $(T_0)_F=gT_Fg^{-1}$.
By \eqref{torus-factorization} for $T_0$, we have $T_0(F)\subset C\cdot G(\mathcal O)$.
Conjugation by $g$ preserves $C\cdot G(\mathcal O)$, and $gS_Fg^{-1}\subset (T_0)_F$.
So $S(F)\subset C\cdot G(\mathcal O)$.

Choose a cocharacter $\lambda\colon\bG_{m,\mathcal O}\to S$ whose weights on $\Lie(U)$ are positive.
For every $u\in U(F)$, the elements $\lambda(\varpi^m)u\lambda(\varpi^m)^{-1}$ tend to $1$ as $m\to\infty$.
For $m\gg 0$ they lie in the open subgroup $N$, and the normality of $N$ forces $u\in N\subset C$.
The same argument with the opposite cocharacter proves $U_Q(F)\subset C$.

For the anisotropic quotient $L/S$, over the valuation rings $\wh{R_i}$ we have $(L/S)(F)=(L/S)(\mathcal O)$.
By the splitting of $S$ and Hilbert's 90 over the semilocal rings $\mathcal O$ and $F$, we have $L(F)=S(F)\cdot L(\mathcal O)$.
Hence $L(F)\subset C\cdot G(\mathcal O)$, and also $P(F)=U(F)L(F)\subset C\cdot G(\mathcal O)$.
Finally, by \cite{SGA3IIInew}*{exposé~XXVI, propositions~4.3.2 and~5.2}, $G(F)\subset U(F)U_Q(F)P(F)$.
Therefore $G(F)\subset C\cdot G(\mathcal O)$.
The reverse inclusion follows from $C\subset G(F)$ and $G(\mathcal O)\subset G(F)$.
\epf

\bpf[Proof of \Cref{decomp-gp}]
Under Setup~\ref{pd-setup}, let $I$ be the image of $G(R_\varpi)$ in $G(F)$, and let $C$ be its closure.
By \Cref{closed-decomp-gp}, $G(F)=C\cdot G(\mathcal O)$.
Since $G(\mathcal O)$ is open in $G(F)$ by \cite{Guo24}*{Lemma~3.5(iii)}, every coset $cG(\mathcal O)$ with $c\in C$ meets $I$.
Thus $C\cdot G(\mathcal O)\subset I\cdot G(\mathcal O)$, and \Cref{closed-decomp-gp} implies $G(F)=I\cdot G(\mathcal O)$.
\epf
We return to the proof of \Cref{G-S over semi-local prufer}.
By the reduction above, assume that $R$ has finite Krull dimension $n$.
We proceed by induction on $n$.
The case $n=0$ is the field case.
Let $\cP$ be a $G$-torsor on $R$ that trivializes over $\Frac R$, and let $\fm_1,\ldots,\fm_r$ be the maximal ideals of $R$.
By prime avoidance, choose $\varpi\in R$ with $V(\varpi)=\{\fm_1,\ldots,\fm_r\}$.
Then $R[\frac{1}{\varpi}]$ is again a semilocal Pr\"ufer domain and has Krull dimension $<n$, so the induction hypothesis implies that $\cP_{R[\frac{1}{\varpi}]}$ is trivial.
For $R_i=R_{\fm_i}$, let $\wh{R_i}$ be the $\varpi$-adic completion of $R_i$.
As in Setup~\ref{pd-setup}, each $\wh{R_i}$ is a rank-one complete valuation ring.
The torsor $\cP_{\wh{R_i}}$ is generically trivial, so it is trivial by \cite{Guo24}*{Theorem~1.3}.
The torsor $\cP$ is therefore trivial after base change to $R[\frac{1}{\varpi}]$ and to $\prod_i\wh{R_i}$.
By \Cref{double cosets}, $\cP$ is represented by a double coset, and \Cref{decomp-gp} identifies this double coset with the trivial one.
Thus $\cP$ is trivial, and the induction is complete.

\end{appendix}


\end{document}